\documentclass[12pt]{article}
\usepackage{amssymb}
\usepackage{amsmath}
\usepackage{amsthm}
\usepackage{color}
\usepackage{comment}
\usepackage{graphicx}
\begin{document}

\def\bw{{\bf w}}
\def\bW{{\bf W}}
\def\by{{\bf y}}
\def\bY{{\bf Y}}
\def\ee{\varepsilon}
\def\myproof{\noindent {\bf Proof.}\\}
\def\re{{\rm e}}
\newcommand{\removableFootnote}[1]{}

\newtheorem{theorem}{Theorem}
\newtheorem{conjecture}[theorem]{Conjecture}
\newtheorem{lemma}[theorem]{Lemma}




\title{The Positive Occupation Time of Brownian Motion with Two-Valued Drift
and Asymptotic Dynamics of Sliding Motion with Noise.}
\author{
D.J.W.~Simpson$^{\dagger}$ and R.~Kuske$^{\ddagger}$\\\\
$^{\dagger}$Institute of Fundamental Sciences\\
Massey University\\
Palmerston North\\
New Zealand\\\\
$^{\ddagger}$Department of Mathematics\\
University of British Columbia\\
Vancouver, BC\\
Canada
}
\maketitle

\begin{abstract}
We derive the probability density function
of the positive occupation time of one-dimensional Brownian motion with two-valued drift.
Long time asymptotics of the density are also computed.
We use the result to describe the transitional probability density function
of a general $N$-dimensional system of stochastic differential equations
representing stochastically perturbed sliding motion of a discontinuous, piecewise-smooth vector field on short time frames.
A description of the density at larger times is obtained via an asymptotic expansion of the Fokker-Planck equation.
\end{abstract}


\section{Introduction}
\label{sec:INTRO}
\setcounter{equation}{0}

Filippov systems are vector fields that are
discontinuous on codimension-one surfaces, termed switching manifolds,
and are used to model a wide variety of physical systems involving a discontinuity or abrupt change, \cite{Fi88,DiBu08}.
Whenever a forward orbit of a Filippov system reaches a section of a switching manifold
for which the vector field points toward the manifold from both sides,
subsequent forward evolution is constrained to the manifold for some time.
Such evolution is called {\em sliding motion}.
In applications, sliding motion often has an important physical interpretation.
For instance, oscillators subject to dry friction
exhibit nonsmooth dynamics when there is both slipping and sticking behaviour \cite{BlCz99,WiDe00,OeHi96,FeMo94}.
Here sliding motion corresponds to the sticking phase of the dynamics.
Models of population dynamics involving assumptions that
species make different selections between habitats or food sources at different times
are often nonsmooth \cite{DeGr07,DeFe06,AmOl12,TaLi12}\removableFootnote{
Also \cite{PiPo13}.
}.
In \cite{DeGr07}, for example, sliding motion is interpreted as the scenario that predators hesitate
between two different food sources.
Also, relay control systems undergo switching
when a controlling signal reaches a threshold value \cite{Ts84,Li03,AsMu08}.
Consequently relay control systems are often well-modelled by piecewise-smooth systems
and sliding motion corresponds to the idealized limit of repeated switching events occurring instantaneously \cite{Jo03,JoRa99,DiJo01}.
In this context sliding motion has been directly exploited to provide enhanced performance \cite{TaLa12,Su06}.

Naturally it is important to understand the effects that noise or uncertainty may have on the dynamics of a Filippov system.
With the addition of small noise in the form of additive Brownian motion, roughly speaking,
sliding motion along a switching manifold becomes a nearby random motion repeatedly intersecting the switching manifold \cite{SiKu13c}.
In general, expressions for the transitional PDF (probability density function)
of a Filippov system with noise are unavailable because the
discontinuity in the drift is an extreme nonlinearity that inhibits exact calculations\removableFootnote{
This is mentioned in \cite{WuHo11}.
}.
In this paper we derive two explicit asymptotic expressions for the transitional PDF of stochastically perturbed sliding motion.
One expression is derived in the context of occupation times and is valid for short time frames;
the other expression is derived via an asymptotic expansion and is valid for time frames that are long
relative to the size of the noise.

As illustrated below in \S\ref{sub:SHORT}, over short time frames we find that
the distance covered by a sample solution in a direction parallel to a switching manifold,
is almost completely determined by the occupation time of the solution on each side of the switching manifold.
Furthermore, over short time frames dynamics in the direction orthogonal to the switching manifold
may be well-approximated by {\em Brownian motion with two-valued drift}:
\begin{equation}
dx(t) =
\left\{ \begin{array}{lc}
a_L \;, & x(t)<0 \\
-a_R \;, & x(t)>0
\end{array} \right\} \,dt + dW(t) \;,
\qquad x(0) = x_0 \;.
\label{eq:dx}
\end{equation}
Here $a_L, a_R \in \mathbb{R}$ and $W(t)$ is a standard Brownian motion.
For any $t>0$, let
\begin{equation}
\tau = \int_0^t \chi_{[0,\infty)}(x(s)) \,ds =
{\rm meas} \left\{ s \in [0,t] ~\big|~ x(s) \ge 0 \right\} \;,
\label{eq:tau}
\end{equation}
denote the {\em positive occupation time} of $x(t)$.
In \S\ref{sub:SHORT} we discuss the approximations more carefully
and use the PDF of $\tau$ to describe the transitional PDF for stochastically perturbed sliding motion
in directions parallel to the switching manifold,
over short time frames.

In this paper we derive the PDF of $\tau$, call it $p$, i.e.,
\begin{equation}
{\rm Prob}(\tau \in \Sigma) = \int_\Sigma p(\hat{\tau};t;x_0,a_L,a_R) \,d\hat{\tau} \;,
\end{equation}
for any measurable subset, $\Sigma \subset [0,t]$.
PDFs of occupation times of simple stochastic processes
have found applications in control problems and mathematical finance \cite{Pe99}.
Equation (\ref{eq:dx}) has also arisen in a stochastic control problem \cite{BeSh80},
and general theoretical settings \cite{GrHe01,QiZh02}.
In the context of sliding motion we require, $a_L, a_R > 0$, such that
solutions to (\ref{eq:dx}) rarely stray far from the origin,
but for generality we study (\ref{eq:dx}) without this restriction.
An explicit expression for the PDF of $x(t)$ was first derived by Karatzas and Shreve in \cite{KaSh84}.

When $x_0 = 0$, $p$ is given by the following theorem.

\begin{theorem}~\\
For $x_0=0$ and any $t>0$, the PDF of the positive occupation time (\ref{eq:tau})
of Brownian motion with two-valued drift (\ref{eq:dx}), is
\begin{eqnarray}
p(\tau;t;0,a_L,a_R) &=&
\frac{\re^{\frac{-a_L^2 (t-\tau)}{2}} \re^{\frac{-a_R^2 \tau}{2}}}
{\pi \sqrt{\tau(t-\tau)}} -
\frac{a_R \,\re^{\frac{-a_L^2 (t-\tau)}{2}} \,{\rm erfc}
\left( \frac{a_R \sqrt{\tau}}{\sqrt{2}} \right)}
{\sqrt{2 \pi (t-\tau)}} -
\frac{a_L \,\re^{\frac{-a_R^2 \tau}{2}} \,{\rm erfc}
\left(\frac{a_L \sqrt{t-\tau}}{\sqrt{2}} \right)}
{\sqrt{2 \pi \tau}} \nonumber \\
&&+~
\frac{\sqrt{2} (a_L+a_R)}
{\sqrt{\pi t}}
\,\re^{\frac{-\left( (a_L+a_R) \tau - a_L t \right)^2}{2t}}
\,{\rm erfc} \left( \frac{-(a_L+a_R) \sqrt{\tau(t-\tau)}}{\sqrt{2t}} \right) \nonumber \\
&&+~
\mathcal{F}(\tau;t;a_L,a_R) + \mathcal{F}(t-\tau;t;a_R,a_L) \;,
\label{eq:p}
\end{eqnarray}
where
\begin{eqnarray}
\mathcal{F}(\tau;t;a_L,a_R) &=&
\frac{a_L (2 a_L + a_R)}{2 \sqrt{\pi}} \int_0^{t-\tau}
-\frac{\sqrt{\tau} \re^{\frac{-a_R^2 \tau}{2}} \re^{\frac{-a_L^2 z}{2}}}
{\sqrt{\pi z} (z+\tau)} \nonumber \\
&&+~\frac{a_L z - a_R \tau}{\sqrt{2} (z+\tau)^\frac{3}{2}}
\,\re^{\frac{-\left( a_L z - a_R \tau \right)^2}{2(z+\tau)}}
\,{\rm erfc} \left( \frac{-(a_L+a_R) \sqrt{z \tau}}{\sqrt{2(z+\tau)}} \right) \,dz \;.
\label{eq:Fcal}
\end{eqnarray}
\label{th:p}
\end{theorem}

In \S\ref{sec:PROOF} we prove Theorem \ref{th:p} by applying the Feynman-Kac formula and taking inverse Laplace transforms.
In \S\ref{sec:BEH} we describe $p$ when $x_0 \ne 0$,
and in \S\ref{sub:SPEC} we discuss the special cases, $a_L = a_R = 0$ and $a_L = -a_R$.
The long time asymptotics of (\ref{eq:p}) are described in \S\ref{sub:ASY}, with a focus on the case $a_L, a_R > 0$.

In \S\ref{sec:SLIDE} we introduce an $N$-dimensional system of stochastic differential equations
that describe stochastically perturbed sliding motion.
For short time frames, in \S\ref{sub:SHORT} we approximate the dynamics
by a one-dimensional stochastic differential equation of the form, (\ref{eq:dx}),
and an algebraic equation involving occupation times for the remaining $N-1$ components of the system,
and apply Theorem \ref{th:p}.
An entirely different methodology is required for long time frames.
In \S\ref{sub:LONG} we use the associated Fokker-Planck equation
to derive the leading order term of an asymptotic expansion of the $N$-dimensional transitional PDF
for stochastically perturbed sliding motion.
This requires matching asymptotics to the third level in the expansion
and imposing a consistency condition at the switching manifold.
The resulting expression is a useful approximation to the PDF
when the noise amplitude is small and at times that are long relative to the magnitude of the noise.
In \S\ref{sub:EXAMPLE} we use Monte-Carlo simulations to illustrate the utility of
both the short time and long time approximations.
Finally \S\ref{sec:CONC} contains concluding remarks.

\section{Proof of Theorem \ref{th:p}}
\label{sec:PROOF}
\setcounter{equation}{0}

For all $t > 0$, the stochastic differential equation (\ref{eq:dx}) has
a unique strong solution \cite{Fl11,PrSh98,StVa69,KrRo05}
and therefore the expectation
\begin{equation}
u(x_0,t,\lambda) = \mathbb{E}_{x_0} \left(
\re^{-\lambda \int_0^t \chi_{[0,\infty)} (x(s)) \,ds} \right) \;,
\label{eq:u}
\end{equation}
where $\lambda > 0$, is well-defined.
Equation (\ref{eq:p}) may be obtained from $u$ because we have
\begin{equation}
u(0,t,\lambda) = \int_0^t \re^{-\lambda \tau} p(\tau;t;0,a_L,a_R) \,d\tau \;.
\label{eq:up}
\end{equation}

Following a method applied to similar problems, see for instance \cite{KaSh84,Ak95,St01},
by the Feynman-Kac formula \cite{KaSh91,Ok03}\removableFootnote{
{\bf The Feynman-Kac formula}:
With certain smoothness and growth constraints on functions $a$, $b$, $f$ and $g$,
for the stochastic differential equation
\begin{equation}
dx(t) = a(x(t)) \,dt + b(x(t)) \,dW(t) \;,
\end{equation}
if
\begin{equation}
u(x_0,t) = \mathbb{E}_{x_0} \left( f(x(t)) \re^{-\int_0^t g(x(s)) \,ds} \right) \;,
\end{equation}
then $u(x_0,t)$ is the unique solution to
\begin{equation}
u_t = a(x_0) u_{x_0} + \frac{b(x_0)^2}{2} u_{x_0 x_0} - g(x_0) u \;,
\end{equation}
with initial condition, $u(x_0,0) = f(x_0)$,
on the domain $\mathbb{R} \times \mathbb{R}^+$.
}\removableFootnote{
A readable text on stochastic calculus is \cite{Ba98}, which I have a copy of.
},
$u$ is the unique, bounded, continuous solution to the initial value problem
\begin{equation}
u_t = \left\{ \begin{array}{lc}
a_L u_{x_0} + \frac{1}{2} u_{x_0 x_0} \;, & x_0 < 0 \\
-a_R u_{x_0} + \frac{1}{2} u_{x_0 x_0} - \lambda u \;, & x_0 > 0
\end{array} \right. \;, \qquad u(x_0,0,\lambda) = 1 \;.
\label{eq:PDE}
\end{equation}
Integration of the PDE with respect to $x_0$ over an arbitrarily small neighbourhood
of zero reveals that $u_{x_0}$ is also continuous.

To solve (\ref{eq:PDE}), we take the Laplace transform
\begin{equation}
U(x_0,\omega,\lambda) = \int_0^\infty \re^{-\omega t} u(x_0,t,\lambda) \,dt \;,
\label{eq:U}
\end{equation}
to produce the piecewise-linear ordinary differential equation
\begin{equation}
-1 = \left\{ \begin{array}{lc}
\frac{1}{2} U_{x_0 x_0} + a_L U_{x_0} - \omega U \;, & x_0 < 0 \\
\frac{1}{2} U_{x_0 x_0} - a_R U_{x_0} - (\omega+\lambda) U \;, & x_0 > 0
\end{array} \right. \;.
\label{eq:ODE}
\end{equation}
Via standard ODE methods, we can obtain an explicit expression
for the unique bounded solution to (\ref{eq:ODE})
for which $u$ and $u_{x_0}$ are continuous at $x_0 = 0$\removableFootnote{
The general solution to the ODE is
\begin{equation}
U(x_0,\omega,\lambda) = \left\{ \begin{array}{lc}
\frac{1}{\omega} +
c_1 \,\re^{x_0 \left( \sqrt{2 \omega + a_L^2} - a_L \right)} +
c_2 \,\re^{-x_0 \left( \sqrt{2 \omega + a_L^2} + a_L \right)} \;, & x_0 < 0 \\
\frac{1}{\omega+\lambda} +
c_3 \,\re^{x_0 \left( \sqrt{2 (\omega+\lambda) + a_R^2} + a_R \right)} +
c_4 \,\re^{-x_0 \left( \sqrt{2 (\omega+\lambda) + a_R^2} - a_R \right)} \;, & x_0 < 0
\end{array} \right. \;,
\end{equation}
for some constants $c_1$, $c_2$, $c_3$ and $c_4$.
For the solution to be bounded we must have $c_2 = c_3 = 0$
(assuming $\omega, \lambda > 0$) regardless of the values of $a_L$ and $a_R$.
To determine $c_1$ and $c_4$, we appeal to the fact that
$u$ and $u_{x_0}$ are continuous at $x_0 = 0$
and obtain the relations
\begin{eqnarray}
\frac{1}{\omega} + c_1 &=& \frac{1}{\omega+\lambda} + c_4 \;, \\
c_1 \left( \sqrt{2 \omega + a_L^2} - a_L \right) &=&
-c_4 \left( \sqrt{2 (\omega+\lambda) + a_R^2} - a_R \right) \;,
\end{eqnarray}
and thus
\begin{eqnarray}
c_1 &=& -\left( \frac{1}{\omega} - \frac{1}{\omega+\lambda} \right)
\frac{\sqrt{2 (\omega+\lambda) + a_R^2} - a_R}
{\sqrt{2 (\omega+\lambda) + a_R^2} - a_R + \sqrt{2 \omega + a_L^2} - a_L} \;, \\
c_4 &=& \left( \frac{1}{\omega} - \frac{1}{\omega+\lambda} \right)
\frac{\sqrt{2 \omega + a_L^2} - a_L}
{\sqrt{2 (\omega+\lambda) + a_R^2} - a_R + \sqrt{2 \omega + a_L^2} - a_L} \;.
\end{eqnarray}
},
and from this arrive at
\begin{equation}
U(0,\omega,\lambda) =
\frac{\omega \left( \sqrt{2 (\omega+\lambda) + a_R^2} - a_R \right) +
(\omega+\lambda) \left( \sqrt{2 \omega + a_L^2} - a_L \right)}
{\omega (\omega+\lambda) \left( \sqrt{2 (\omega+\lambda) + a_R^2} - a_R +
\sqrt{2 \omega + a_L^2} - a_L \right)} \;.
\label{eq:U0}
\end{equation}
Our goal is to obtain $p$,
which by (\ref{eq:up}) and (\ref{eq:U}) is related to $U(0,\omega,\lambda)$ by
\begin{equation}
U(0,\omega,\lambda) = \int_0^\infty \int_0^t
\re^{-\omega t} \re^{-\lambda \tau}
p(\tau;t;0,a_L,a_R) \,d\tau \,dt \;.
\label{eq:Up}
\end{equation}
Equation (\ref{eq:Up}) is now used to derive $p$ {\em constructively}.
To make the task simply the evaluation of two inverse Laplace transforms, 
we first reverse the order of integration in (\ref{eq:Up}),
and let $v = t-\tau$ and $\mu = \omega + \lambda$, to obtain
\begin{equation}
U(0,\omega,\mu-\omega) = \int_0^\infty \int_0^\infty
\re^{-\omega v} \re^{-\mu \tau}
p(\tau,v+\tau;0,a_L,a_R) \,dv \,d\tau \;.
\label{eq:Up2}
\end{equation}
Simplification is provided by noting from (\ref{eq:U0}) that we can write
\begin{equation}
U(0,\omega,\mu-\omega) = \Xi(\mu,\omega;a_L,a_R) + \Xi(\omega,\mu;a_R,a_L) \;,
\end{equation}
where
\begin{equation}
\Xi(\mu,\omega;a_L,a_R) =
\frac{\sqrt{2 \mu + a_R^2} - a_R}
{\mu \left( \sqrt{2 \mu + a_R^2} - a_R + \sqrt{2 \omega + a_L^2} - a_L \right)} \;.
\label{eq:Xi}
\end{equation}
By (\ref{eq:Up2}), with a little care it follows that
\begin{equation}
p(\tau;t;0,a_L,a_R) = \xi(\tau,t-\tau;a_L,a_R) + \xi(t-\tau,\tau;a_R,a_L) \;,
\label{eq:p2}
\end{equation}
where
\begin{equation}
\Xi(\mu,\omega;a_L,a_R) =
\int_0^\infty \int_0^\infty \re^{-\omega v} \re^{-\mu \tau}
\xi(\tau,v;a_L,a_R) \,dv \,d\tau \;.
\label{eq:Xi2}
\end{equation}
A determination of $\xi$ from (\ref{eq:Xi2}) using (\ref{eq:Xi})
is achieved by inverting the double Laplace transform.
The details of this calculation are given in Appendix \ref{sec:INVLAP}.
The result for $\xi$ combined with (\ref{eq:p2})
then produces (\ref{eq:p}) with (\ref{eq:Fcal}) as required.
\hfill
$\Box$

\section{The behaviour and asymptotics of $p$}
\label{sec:BEH}
\setcounter{equation}{0}

Theorem \ref{th:p} gives the PDF, $p$, of the positive occupation time, $\tau$, when $x_0 = 0$ (where $x_0 = x(0)$).
The result for $x_0 \ne 0$ may be written in terms of the result for $x_0 = 0$ by
applying the strong Markov property of (\ref{eq:dx}).
Specifically, let
\begin{equation}
h(s;x_0,a) = \frac{|x_0|}{\sqrt{2 \pi s^3}}
\,\re^{\frac{-(x_0+as)^2}{2 s}} \;,
\label{eq:h}
\end{equation}
denote the PDF for the first passage time, $s$, of
$dx(t) = a \,dt + dW(t)$,
from $x(0)=x_0$ to $x(t)=0$.
Then
\begin{equation}
p(\tau;t;x_0,a_L,a_R) =
\left\{ \begin{array}{lc}
\int_0^t h(s;x_0,a_L) p(\tau;t-s;0,a_L,a_R) \,ds +
\delta(\tau) \int_t^\infty h(s,x_0,a_L) \,ds \;, & x_0 < 0 \\
\int_0^t h(s;x_0,-a_R) p(\tau-s;t-s;0,a_L,a_R) \,ds \\
\qquad+~\delta(\tau-t) \int_t^\infty h(s,x_0,-a_R) \,ds \;, & x_0 > 0
\end{array} \right. \;,
\label{eq:px0nonzero}
\end{equation}
where $\delta$ is the Dirac-delta function
and we set $p(\tau;t;0,a_L,a_R) = 0$ whenever $\tau \notin [0,t]$.
For both $x_0 < 0$ and $x_0 > 0$,
the first term of (\ref{eq:px0nonzero}) corresponds to the case
that $x(t)$ first reaches zero at a time $s \le t$,
and the second term corresponds to $x(t)$ not reaching zero by the time $t$.

\subsection{Two special cases}
\label{sub:SPEC}

\begin{figure}[t!]
\begin{center}
\setlength{\unitlength}{1cm}
\begin{picture}(16,8.5)
\put(0,.2){\includegraphics[height=8cm]{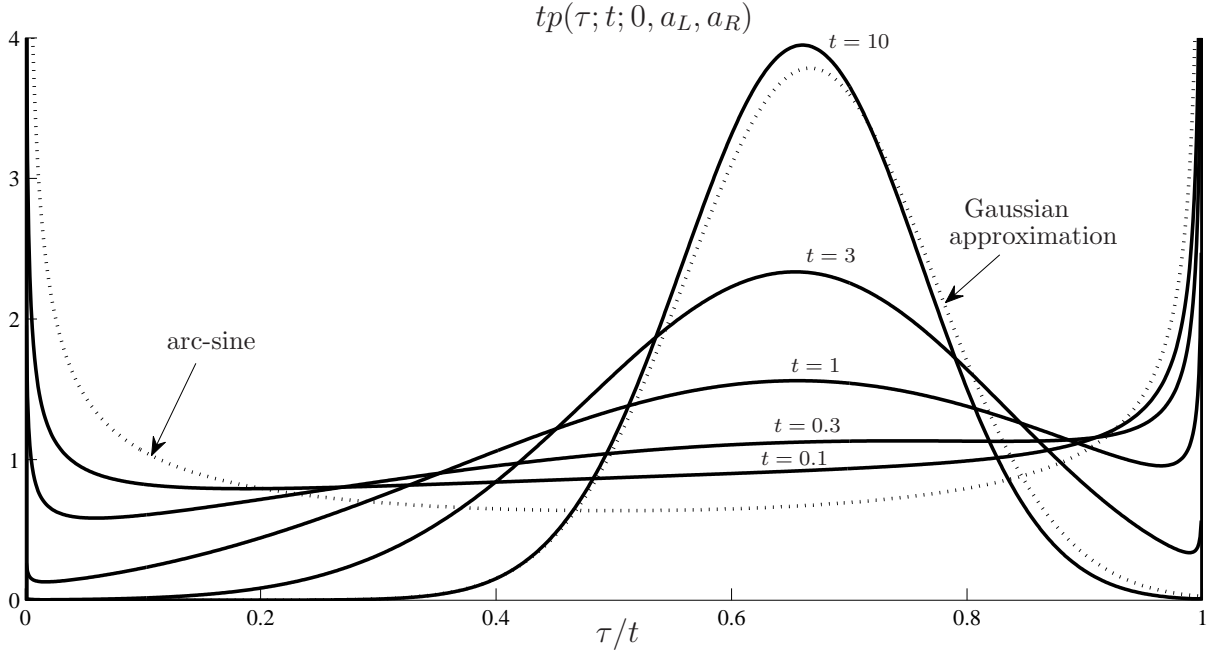}}
\put(7.8,0){$\tau/t$}
\put(7,8.2){$t p(\tau;t;0,a_L,a_R)$}
\put(10,2.36){\scriptsize $t = 0.1$}
\put(10.17,2.8){\scriptsize $t = 0.3$}
\put(10.4,3.6){\scriptsize $t = 1$}
\put(10.6,5.05){\scriptsize $t = 3$}
\put(10.9,7.9){\scriptsize $t = 10$}
\put(2.1,3.9){\footnotesize arc-sine}
\put(12.7,5.65){\footnotesize Gaussian}
\put(12.5,5.3){\footnotesize approximation}
\end{picture}
\caption{
The PDF, $p$, when $a_L = 2$ and $a_R = 1$,
computed by numerically evaluating the integrals in (\ref{eq:p})-(\ref{eq:Fcal}).
Note that the axes are scaled by $t$ such that the domain
is the unit interval at all times.
For small $t$, $p$ is well-approximated by the density
of the arc-sine distribution, (\ref{eq:Le39}).
For large $t$, $p$ is approximately
the Gaussian density of (\ref{eq:pLarget}), since $a_L, a_R > 0$.
When $t = 10$, on the given axes the Gaussian approximation has mean
$\frac{a_L}{a_L+a_R} = \frac{2}{3}$
and standard deviation
$\frac{1}{(a_L+a_R) \sqrt{t}} = \frac{1}{3 \sqrt{10}} \approx 0.11$,
as shown.
\label{fig:occTimePDF}
}
\end{center}
\end{figure}

In two special cases, (\ref{eq:p}) reduces to previously described results.
First, if $a_L = a_R = 0$, then $x(t)$ is simply regular Brownian motion.
As first shown by Levy \cite{Le39},
in this case the occupation time is distributed by the
arc-sine distribution \cite{KaSh91,Bi99}.
Indeed, in this case only the first term of (\ref{eq:p}) is nonzero, and
\begin{equation}
p(\tau;t;0,0,0) = \frac{1}{\pi \sqrt{\tau(t-\tau)}} \;.
\label{eq:Le39}
\end{equation}
Equation (\ref{eq:Le39}) provides a reasonable approximation
to $p(\tau;t;0,a_L,a_R)$ for small $t$, see Fig.~\ref{fig:occTimePDF},
because over short time intervals the diffusion of a stochastic quantity dominates its drift.

Second, if $a_L = -a_R$, then $x(t)$ is Brownian motion with constant drift.
In this case (\ref{eq:p}) simplifies to
\begin{equation}
p(\tau;t;0,-a,a) =
\left( \frac{1}{\sqrt{\pi \tau}} \,\re^{\frac{-a^2 \tau}{2}} -
\frac{a}{\sqrt{2}} \,{\rm erfc} \left( \frac{a \sqrt{\tau}}{\sqrt{2}} \right) \right)
\left( \frac{1}{\sqrt{\pi (t-\tau)}} \,\re^{\frac{-a^2 (t-\tau)}{2}} +
\frac{a}{\sqrt{2}} \,{\rm erfc} \left( \frac{-a \sqrt{t-\tau}}{\sqrt{2}} \right) \right) \;,
\label{eq:Ak95}
\end{equation}
which was first derived by Akahori \cite{Ak95}.
The simplification is not straight-forward
but may be demonstrated by using different expressions for
Owen's T-function \cite{Ow56}\removableFootnote{
With $-a_L = a_R = a$, the first three terms of (\ref{eq:p}) match three
of the four terms of (\ref{eq:Ak95}) that result from expanding brackets.
The fourth term of (\ref{eq:p}) vanishes, thus it remains to show that
\begin{eqnarray}
\mathcal{F}(\tau;t;-a,a) + \mathcal{F}(t-\tau;t;a,-a) &=&
\frac{a^2}{2 \sqrt{\pi}} \int_0^{t-\tau}
\left( -\frac{a}{\sqrt{2}} - \frac{\sqrt{\tau}}{\sqrt{\pi z} \sqrt{z+\tau}} \right)
\frac{\re^{\frac{-a^2(z+\tau)}{2}}}{\sqrt{z+\tau}} \,dz \nonumber \\
&&+~
\frac{a^2}{2 \sqrt{\pi}} \int_0^\tau
\left( \frac{a}{\sqrt{2}} - \frac{\sqrt{t-\tau}}{\sqrt{\pi z} \sqrt{z+t-\tau}} \right)
\frac{\re^{\frac{-a^2(z+t-\tau)}{2}}}{\sqrt{z+t-\tau}} \,dz \nonumber \\
&=&
-a^2 {\rm erfc} \left( \frac{a \sqrt{\tau}}{\sqrt{2}} \right) -
\frac{a^2 \sqrt{\tau}}{2 \pi} \int_0^{t-\tau}
\frac{\re^{\frac{-a^2 (z+\tau)}{2}}}{\sqrt{z}(z+\tau)} \,dz \nonumber \\ 
&&+~
a^2 {\rm erfc} \left( \frac{a \sqrt{t-\tau}}{\sqrt{2}} \right) -
\frac{a^2 \sqrt{t-\tau}}{2 \pi} \int_0^\tau
\frac{\re^{\frac{-a^2 (z+t-\tau)}{2}}}{\sqrt{z}(z+t-\tau)} \,dz \;,
\label{eq:constantDriftLastTerm}
\end{eqnarray}
simplifies to the remaining term of (\ref{eq:Ak95}):
$-\frac{a^2}{2} \,{\rm erfc} \left( \frac{a \sqrt{\tau}}{\sqrt{2}} \right)
\,{\rm erfc} \left( \frac{-a \sqrt{t-\tau}}{\sqrt{2}} \right)$.
This fact may be demonstrated by rewriting the integrals as double
integrals of the density of a bivariate Gaussian
(Owen's T-function is
$T(h,\alpha) = \frac{1}{2 \pi} \int_0^\alpha
\frac{\re^{\frac{-h^2 (1+x^2)}{2}}}{1+x^2} \,dx =
\int_h^\infty \int_0^{\alpha x}
\frac{1}{2 \pi} \,\re^{\frac{-(x^2+y^2)}{2}} \,dy \,dx$),
or differentiating and integrating with respect to $b \equiv a^2$:
\begin{eqnarray}
&& \frac{\partial}{\partial b} \Bigg(
\frac{2}{b} \bigg( \mathcal{F} \left( \tau;t;-\sqrt{b},\sqrt{b} \right) 
+ \mathcal{F} \left( t-\tau;t;\sqrt{b},-\sqrt{b} \right) \bigg) \Bigg) \nonumber \\
&&=~
\frac{\sqrt{\tau}}{\sqrt{\pi b}} \re^{\frac{-b \tau}{2}} +
\frac{\sqrt{\tau}}{\pi} \int_0^{t-\tau}
\frac{\re^{\frac{-b(z+\tau)}{2}}}{\sqrt{z}} \,dz -
\frac{\sqrt{\tau}}{\sqrt{\pi b}} \re^{\frac{-b (t-\tau)}{2}} +
\frac{\sqrt{\tau}}{\pi} \int_0^\tau
\frac{\re^{\frac{-b(z+t-\tau)}{2}}}{\sqrt{z}} \,dz \nonumber \\
&&=~
\frac{\sqrt{\tau}}{\sqrt{\pi b}} \re^{\frac{-b \tau}{2}}
{\rm erfc} \left( \frac{-\sqrt{b (t-\tau)}}{\sqrt{2}} \right) -
\frac{\sqrt{t-\tau}}{\sqrt{\pi b}} \re^{\frac{-b (t-\tau)}{2}}
{\rm erfc} \left( \frac{\sqrt{b \tau}}{\sqrt{2}} \right) \nonumber \\
&&=~
\frac{\partial}{\partial b} \left( -
{\rm erfc} \left( \frac{\sqrt{b \tau}}{\sqrt{2}} \right)
{\rm erfc} \left( \frac{-\sqrt{b (t-\tau)}}{\sqrt{2}} \right) \right) \;.
\end{eqnarray}
}.

\subsection{Long time asymptotics of $p$}
\label{sub:ASY}

Here we provide a statement of the long time asymptotics of (\ref{eq:p}), then discuss its derivation.

For large $t$,
\begin{equation}
p(\tau;t;0,a_L,a_R) \sim \left\{ \begin{array}{lc}
\frac{a_L+a_R}{\sqrt{2 \pi t}}
\,\re^{\frac{-\left( (a_L+a_R) \tau - a_L t \right)^2}{2 t}}
\;, & a_L > 0, \,a_R > 0 \\
\mathcal{G}(\tau;a_L,a_R)
\;, & a_L < 0, \,a_R > 0 \\
\mathcal{G}(t-\tau;a_R,a_L)
\;, & a_L > 0, \,a_R < 0 \\
\mathcal{G}(\tau;a_L,a_R) + \mathcal{G}(t-\tau;a_R,a_L)
\;, & a_L < 0, \,a_R < 0
\end{array} \right. \;,
\label{eq:pLarget}
\end{equation}
where,
\begin{equation}
\mathcal{G}(\tau;a_L,a_R) =
-\frac{\sqrt{2} a_L}{\sqrt{\pi \tau}}
\,\re^{\frac{-a_R^2 \tau}{2}} -
a_L (2 a_L + a_R)
\,\re^{2 a_L (a_L+a_R) \tau}
\,{\rm erfc} \left( \frac{-(2 a_L + a_R) \sqrt{\tau}}{\sqrt{2}} \right) \;.
\label{eq:Gcal}
\end{equation}
To obtain (\ref{eq:pLarget}),
note that it may be shown directly from (\ref{eq:p}),
that for any fixed $\tau > 0$,
\begin{equation}
p(\tau;t;0,a_L,a_R) \sim
\left\{ \begin{array}{lc}
\mathcal{G}(\tau;a_L,a_R) \;, & a_L < 0 \\
0 \;, & a_L > 0
\end{array} \right. \;,
\end{equation}
for large $t$.
For brevity we omit a derivation of this statement which
may be demonstrated by
taking $t \to \infty$ in (\ref{eq:Fcal})
and simplifying the result via integral transformations\removableFootnote{
When $a_L < 0$ and $a_R > 0$, as $t \to \infty$ the first, second and fourth
terms of (\ref{eq:p}) vanish, as does $\mathcal{F}(t-\tau;\tau;a_R,a_L)$.
The third term tends to 
$-\frac{\sqrt{2} a_L}{\sqrt{\pi \tau}}
\,\re^{\frac{-a_R^2 \tau}{2}}$.
Thus in order to obtain (\ref{eq:pLarget}) for $a_L < 0$, $a_R > 0$, it remains to show:
\begin{equation}
\lim_{t \to \infty} \mathcal{F}(\tau;t;a_L,a_R) =
-a_L (2 a_L + a_R)
\,\re^{2 a_L (a_L+a_R) \tau}
\,{\rm erfc} \left( \frac{-(2 a_L + a_R) \sqrt{\tau}}{\sqrt{2}} \right) \;.
\end{equation}
To do this we split the integral in $\mathcal{F}$:
\begin{equation}
\frac{2 \sqrt{\pi}}{a_L (2 a_L + a_R)}
\lim_{t \to \infty} \mathcal{F}(\tau;t;a_L,a_R) =
\int_0^\infty
-\frac{\sqrt{\tau} \re^{\frac{-a_R^2 \tau}{2}} \re^{\frac{-a_L^2 z}{2}}}
{\sqrt{\pi z} (z+\tau)} +
\frac{a_L z - a_R \tau}{\sqrt{2} (z+\tau)^\frac{3}{2}}
\,\re^{\frac{-\left( a_L z - a_R \tau \right)^2}{2(z+\tau)}}
\,{\rm erfc} \left( \frac{-(a_L+a_R) \sqrt{z \tau}}{\sqrt{2(z+\tau)}} \right) \,dz \;,
\end{equation}
into pieces in two different ways.
In the subcase $-a_R \le a_L < 0$,
we use ${\rm erfc}(-X) = 2 - {\rm erfc}(X)$ to write
\begin{eqnarray}
\frac{2 \sqrt{\pi}}{a_L (2 a_L + a_R)}
\lim_{t \to \infty} \mathcal{F}(\tau;t;a_L,a_R) 
&=& -\frac{\sqrt{\tau} \re^{\frac{-a_R^2 \tau}{2}}}{\sqrt{\pi}} \int_0^\infty
\frac{\re^{\frac{-a_L^2 z}{2}}}{\sqrt{z}(z+\tau)} \,dz \nonumber \\
&&+~
\sqrt{2} \int_0^\infty
\frac{a_L z - a_R \tau}{(z+\tau)^\frac{3}{2}}
\,\re^{\frac{-\left( a_L z - a_R \tau \right)^2}{2(z+\tau)}} \,dz \nonumber \\
&&-~
\frac{1}{\sqrt{2}} \int_0^\infty
\frac{a_L (z+\tau) - (a_L+a_R) \tau}{(z+\tau)^\frac{3}{2}}
\,\re^{\frac{-\left( a_L z - a_R \tau \right)^2}{2(z+\tau)}}
\,{\rm erfc} \left( \frac{(a_L+a_R) \sqrt{z \tau}}{\sqrt{2(z+\tau)}} \right) \,dz \;,
\end{eqnarray}
and if $0 < a_R \le -a_L$, we use 
$a_L z - a_R \tau = -a_L (z+\tau) - (a_L+a_R) \tau + 2 a_L (z+\tau)$ to write
\begin{eqnarray}
\frac{2 \sqrt{\pi}}{a_L (2 a_L + a_R)}
\lim_{t \to \infty} \mathcal{F}(\tau;t;a_L,a_R) 
&=& -\frac{\sqrt{\tau} \re^{\frac{-a_R^2 \tau}{2}}}{\sqrt{\pi}} \int_0^\infty
\frac{\re^{\frac{-a_L^2 z}{2}}}{\sqrt{z}(z+\tau)} \,dz \nonumber \\
&&+~
\sqrt{2} a_L \int_0^\infty
\frac{1}{\sqrt{z+\tau}}
\,\re^{\frac{-\left( a_L z - a_R \tau \right)^2}{2(z+\tau)}}
\,{\rm erfc} \left( \frac{-(a_L+a_R) \sqrt{z \tau}}{\sqrt{2(z+\tau)}} \right) \,dz \nonumber \\
&&-~
\frac{1}{\sqrt{2}} \int_0^\infty
\frac{a_L (z+\tau) + (a_L+a_R) \tau}{(z+\tau)^\frac{3}{2}}
\,\re^{\frac{-\left( a_L z - a_R \tau \right)^2}{2(z+\tau)}}
\,{\rm erfc} \left( \frac{-(a_L+a_R) \sqrt{z \tau}}{\sqrt{2(z+\tau)}} \right) \,dz \;.
\end{eqnarray}
The result then follows from Lemma \ref{le:evalIntInf}.
\begin{lemma}
Suppose $\tau > 0$ and $a_L < 0$.
Then
\begin{eqnarray}
\frac{\sqrt{\tau} \re^{\frac{-a_R^2 \tau}{2}}}{\sqrt{\pi}}
\int_0^\infty \frac{\re^{\frac{-a_L^2 z}{2}}}{\sqrt{z}(z+\tau)} \,dz &=& 
\sqrt{\pi} \,\re^{\frac{a_L^2 \tau}{2}} \,\re^{\frac{-a_R^2 \tau}{2}}
\,{\rm erfc} \left( \frac{-a_L \sqrt{\tau}}{\sqrt{2}} \right) \;,
\label{eq:evalIntInf1} \\
\sqrt{2} \int_0^\infty
\frac{a_L z - a_R \tau}{(z+\tau)^{\frac{3}{2}}}
\,\re^{\frac{-(a_L z - a_R \tau)^2}{2(z+\tau)}} \,dz &=&
-2 \sqrt{\pi} \,\re^{2 a_L (a_L+a_R) \tau}
\,{\rm erfc} \left( \frac{-(2 a_L + a_R) \sqrt{\tau}}{\sqrt{2}} \right) \;,
\label{eq:evalIntInf2} 
\end{eqnarray}
and if $a_R \ge -a_L$, then
\begin{equation}
\hspace{-15mm}
\frac{1}{\sqrt{2}} \int_0^\infty
\frac{a_L (z+\tau) - (a_L+a_R) \tau}{(z+\tau)^\frac{3}{2}}
\,\re^{\frac{-\left( a_L z - a_R \tau \right)^2}{2(z+\tau)}}
\,{\rm erfc} \left( \frac{(a_L+a_R) \sqrt{z \tau}}{\sqrt{2(z+\tau)}} \right) \,dz =
-\sqrt{\pi} \,\re^{\frac{a_L^2 \tau}{2}} \,\re^{\frac{-a_R^2 \tau}{2}}
\,{\rm erfc} \left( \frac{-a_L \sqrt{\tau}}{\sqrt{2}} \right) \;,
\label{eq:evalIntInf3}
\end{equation}
and if $a_R \le -a_L$, then
\begin{equation}
\hspace{-15mm}
\sqrt{2} a_L \int_0^\infty
\frac{1}{\sqrt{z+\tau}}
\,\re^{\frac{-\left( a_L z - a_R \tau \right)^2}{2(z+\tau)}}
\,{\rm erfc} \left( \frac{-(a_L+a_R) \sqrt{z \tau}}{\sqrt{2(z+\tau)}} \right) \,dz =
-2 \sqrt{\pi} \,\re^{2 a_L (a_L+a_R) \tau}
\,{\rm erfc} \left( \frac{-(2 a_L + a_R) \sqrt{\tau}}{\sqrt{2}} \right) \;,
\label{eq:evalIntInf4}
\end{equation}
(although numerics indicate that (\ref{eq:evalIntInf3}) and (\ref{eq:evalIntInf4})
hold for any $a_R \in \mathbb{R}$, but my proof requires the given restrictions).
\label{le:evalIntInf}
\end{lemma}

\noindent
{\bf Proof}.~~Equation (\ref{eq:evalIntInf1})
follows from the identity:
$\sqrt{\tau} \,\re^{\frac{-a_L^2 \tau}{2}}
\int_0^v \frac{\re^{\frac{-a_L^2 z}{2}}}{\sqrt{z}(z+\tau)} \,dz =
4 \int_{\frac{|a_L| \sqrt{\tau}}{\sqrt{2}}}^\infty
\int_0^{\frac{\sqrt{v} x}{\sqrt{\tau}}}
\re^{-x^2} \,\re^{-y^2} \,dy \,dx$,
which is essentially Owen's T-function
and may be verified by differentiating both sides with respect to $a_L^2$.

Equation (\ref{eq:evalIntInf2}) may be demonstrated by
splitting the integral in two via
$a_L z - a_R \tau = a_L z + (2 a_L + a_R) \tau - 2 (a_L + a_R) \tau$,
and substituting
$x = \frac{a_L z - a_R \tau}{\sqrt{2(z+\tau)}}$ to the $a_L z + (2 a_L + a_R) \tau$ piece,
and $y = \frac{\tau}{\sqrt{z+\tau}}$ to the $2 (a_L + a_R) \tau$ piece, and utilizing
$\int \re^{-a^2 y^2} \,\re^{\frac{-b^2}{y^2}} \,dy =
-\frac{\sqrt{\pi}}{4 a} \left(
\re^{-2ab} \,{\rm erfc} \left( ay - \frac{b}{y} \right) +
\re^{2ab} \,{\rm erfc} \left( ay + \frac{b}{y} \right) \right) + C$.

When $a_L+a_R \ge 0$, we can write the left hand side of (\ref{eq:evalIntInf3}) as
\begin{equation}
\frac{\re^{\frac{-a_R^2 \tau}{2}}}{\sqrt{\pi}}
\int_0^\infty \int_0^\infty
\frac{a_L(z+\tau)-a_S \tau}{(z+\tau) \sqrt{2(z+\tau)w + a_S^2 z \tau}}
\,\re^{\frac{-a_L^2 z}{2}} \,\re^{-w} \,dw \,dz =
\frac{\re^{\frac{-a_R^2 \tau}{2}}}{\sqrt{\pi}} \mathcal{H}(\tau;a_L,a_S) \;,
\end{equation}
where we define $a_S = a_L+a_R$, and the function $\mathcal{H}$ as indicated.
Equation (\ref{eq:evalIntInf3}) is then easily verified in the case $a_S = 0$.
We can further prove (\ref{eq:evalIntInf3}) holds for any $a_S \ge 0$
by showing that the derivative
\begin{equation}
\frac{\partial \mathcal{H}}{\partial a_S} =
-\frac{\sqrt{\tau}}{\sqrt{2}} \int_0^\infty \int_0^\infty
\frac{\left( \frac{a_S \sqrt{\tau}}{\sqrt{2}} u + \frac{a_L \sqrt{\tau}}{\sqrt{2}} w \right)
\,\re^{-u} \,\re^{-w}}
{\left( \frac{a_S^2 \tau}{2} u + \frac{a_L^2 \tau}{2} w + u w \right)^{\frac{3}{2}}}
\,du \,dw \;,
\label{eq:dHdaS}
\end{equation}
is always zero.
Indeed $\frac{\partial \mathcal{H}}{\partial a_S} = 0$ when $a_L < 0$,
as may be demonstrated by performing the two-dimensional substitution $(u,w) = (x+y,x-y)$,
and showing that the resulting integral in $y$ is identically zero.

To prove (\ref{eq:evalIntInf4})
we subtract it from (\ref{eq:evalIntInf2})
and apply the substitutions 
$x = \frac{(a_L+a_R)^2 \tau^2}{2 (z+\tau)}$,
$b = \frac{(a_L+a_R)^2 \tau}{2}$ and
$c = \frac{a_L (a_L+a_R) \tau}{2}$ to obtain
\begin{equation}
\sqrt{2} \int_0^\infty
\left(
\frac{a_L z - a_R \tau}{(z+\tau)^{\frac{3}{2}}} -
\frac{a_L}{\sqrt{z+\tau}}
\,{\rm erfc} \left( \frac{-(a_L+a_R) \sqrt{z \tau}}{\sqrt{2(z+\tau)}} \right)
\right)
\,\re^{\frac{-(a_L z - a_R \tau)^2}{2(z+\tau)}} \,dz =
2 \re^{2c} \int_0^b
\frac{c - x - c \,{\rm erfc} \left( \sqrt{b-x} \right)}{x^{\frac{3}{2}}}
\,\re^{-x} \,\re^{\frac{-c^2}{x}} \,dx \;.
\end{equation}
The right hand side of this equation is zero when $b = 0$.
Moreover the right hand side is zero for all $b \ge 0$ because
(as is straight-forward to demonstrate)
its derivative with respect to $b$ is identically zero.
}.
Moreover, for $a_L < 0$\removableFootnote{
A demonstration of this statement requires considering
different cases of $a_L$ and $a_R$ separately and carefully.
There results a double integral similar to (\ref{eq:dHdaS})
that may be evaluated via the same substitution used above.
},
\begin{equation}
\int_0^\infty \mathcal{G}(\tau;a_L,a_R) \,d\tau =
\left\{ \begin{array}{lc}
\frac{a_L}{a_L+a_R}  \;, & a_R \le 0 \\
1 \;, & a_R \ge 0
\end{array} \right. \;,
\end{equation}
which ensures (\ref{eq:pLarget}) has unit area in the limit $t \to \infty$.
The term $\mathcal{G}(t-\tau;a_R,a_L)$ appears in (\ref{eq:pLarget}) by symmetry.

When $a_L, a_R > 0$, for large $t$
the first three terms of (\ref{eq:p}) may be neglected 
because their $L_1$ norms rapidly approach zero as $t \to \infty$.
The fourth term of (\ref{eq:p}) approaches
$\frac{4 (a_L+a_R)}{\sqrt{2 \pi t}}
\re^{\frac{-\left( (a_L+a_R) \tau - a_L t \right)^2}{2 t}}$,
because ${\rm erfc}(\cdot) \to 2$\removableFootnote{
Since $\tau = O(t)$ if the exponential term is to be non-negligible.
}.
To simplify (\ref{eq:Fcal}),
we consider the two terms in the integrand of (\ref{eq:Fcal}) separately.
The first term may be neglected for large $t$ when $a_L, a_R > 0$.
To deal with the second term we let
\begin{equation}
\rho = t^\eta \left( \tau - \frac{a_L t}{a_L+a_R} \right) \;,
\label{eq:rho}
\end{equation}
be a scaled difference of $\tau$ from the mean, $\frac{a_L t}{a_L+a_R}$,
for some $\eta \in \mathbb{R}$\removableFootnote{
This must satisfy $\eta \ge -1$ since $\tau \in [0,t]$.
}.
Then, substituting $u = \frac{a_L z - a_R \tau}{t}$ into (\ref{eq:Fcal}) produces
\begin{equation}
\mathcal{F}(\tau;t;a_L,a_R) \sim
\frac{a_L^{\frac{3}{2}} (2 a_L + a_R) \sqrt{t}}{2\sqrt{2 \pi}}
\int_{-\frac{a_L a_R}{a_L+a_R} - \frac{a_R \rho}{t^{\eta+1}}}^{-\frac{(a_L+a_R) \rho}{t^{\eta+1}}}
\frac{u}{\left( a_L + \frac{(a_L+a_R) \rho}{t^{\eta+1}} + u \right)^{\frac{3}{2}}}
\,\re^{\frac{-a_L t u^2}{2 \left( a_L + \frac{(a_L+a_R) \rho}{t^{\eta+1}} + u \right)}}
\,{\rm erfc}(\cdot) \,du \;,
\label{eq:FcalPiece2}
\end{equation}
where the argument of the complement error function is omitted for brevity.
Equation (\ref{eq:FcalPiece2}) may be evaluated asymptotically \cite{BeOr99}\removableFootnote{
Technically this is not Laplace's method \cite{BeOr99}
because we do not take the upper limit of integration to $\infty$.
}.
When $\rho \le 0$, this is achieved by expanding the integrand in a Taylor series about $u=0$;
when $\rho > 0$, we expand about the upper limit of the integral.
For instance, when $\rho \le 0$
(the other case is similar and produces the same result\removableFootnote{
When $\rho > 0$, we can let $v = u + \frac{(a_L+a_R) \rho}{t^{\eta+1}}$,
expand about $v = 0$, and, at least when $-1 \le \eta < 0$
(otherwise the $v$-term in the exponential
is not large and I think we need the $v^2$ term),
the asymptotic evaluation of the integral gives the same result as when $\rho \le 0$.
}),
after integration we have\removableFootnote{
The intermediate step is:
\begin{equation}
\mathcal{F}(\tau;t;a_L,a_R) \sim
\frac{a_L^{\frac{3}{2}} (2 a_L + a_R) \sqrt{t}}
{\sqrt{2 \pi} \left( a_L + \frac{(a_L+a_R) \rho}{t^{\eta+1}} \right)^{\frac{3}{2}}}
\int_{-\frac{a_L a_R}{a_L+a_R} - \frac{a_R \rho}{t^{\eta+1}}}^{-\frac{(a_L+a_R) \rho}{t^{\eta+1}}}
u \,\re^{\frac{-a_L t u^2}{2 \left( a_L + \frac{(a_L+a_R) \rho}{t^{\eta+1}} \right)}}
\left( 1 + O(u) \right) \,du \;.
\label{eq:FcalPiece2b}
\end{equation}
}
\begin{equation}
\mathcal{F}(\tau;t;a_L,a_R) \sim
-\frac{\sqrt{a_L} (2 a_L + a_R)}
{\sqrt{2 \pi \left( a_L + \frac{(a_L+a_R) \rho}{t^{\eta+1}} \right) t}}
\,\re^{\frac{-a_L (a_L+a_R)^2 \rho^2}
{2 t^{2 \eta + 1} \left( a_L + \frac{(a_L+a_R) \rho}{t^{\eta+1}} \right)}}
+ O \left( \frac{1}{t} \right) \;.
\label{eq:FcalPiece2c}
\end{equation}
Only when $\eta = -\frac{1}{2}$ is the leading order term of (\ref{eq:FcalPiece2c})
non-constant in the limit $t \to \infty$,
and in this case, after also substituting back $\tau$ using (\ref{eq:rho}),
\begin{equation}
\mathcal{F}(\tau;t;a_L,a_R) \sim
-\frac{2 a_L + a_R}{\sqrt{2 \pi t}}
\,\re^{\frac{-\left( (a_L+a_R) \tau - a_L t \right)^2}{2 t}}
+ O \left( \frac{1}{t} \right) \;.
\label{eq:FcalPiece2d}
\end{equation}
In view of (\ref{eq:FcalPiece2d}),
summing the asymptotic limits of the fourth, fifth and sixth terms of (\ref{eq:p})
yields (\ref{eq:pLarget}) for $a_L, a_R > 0$.

When $a_L, a_R > 0$, (\ref{eq:pLarget}) is the PDF of a Gaussian
with mean $\frac{a_L t}{a_L + a_R}$,
and standard deviation $\frac{\sqrt{t}}{a_L+a_R}$, see Fig.~\ref{fig:occTimePDF}.
In this case the direction of the drift is toward $x=0$
for both positive and negative $x$.
The fraction of time spent in $[0,\infty)$, approaches $\frac{a_L}{a_L+a_R}$ as $t \to \infty$.
This is consistent with the fact that the area under the steady-state density
of $x(t)$ over $[0,\infty)$ is $\frac{a_L}{a_L+a_R}$ \cite{SiKu13c}.

\section{An application to stochastically perturbed sliding motion}
\label{sec:SLIDE}
\setcounter{equation}{0}

Our interest is in the dynamics near a switching manifold of a general $N$-dimensional Filippov system.
Assuming the manifold is smooth, we may choose our coordinate system such that the manifold
coincides with points in $\mathbb{R}^N$ that are zero in their first coordinate \cite{DiBu08,DiBu01}.
We let $x \in \mathbb{R}$ denote the first coordinate and $\by \in \mathbb{R}^{N-1}$ denote the remaining coordinates.
If we ignore any other switching manifolds,
and add constant noise to the system in the form of additive Brownian motion, then the system may be written as
\begin{equation}
\left[ \begin{array}{c} dx(t) \\ d\by(t) \end{array} \right] =
\left\{ \begin{array}{lc}
\left[ \begin{array}{c} \phi^{(L)}(x(t),\by(t)) \\ \psi^{(L)}(x(t),\by(t)) \end{array} \right] \;, & x(t) < 0 \\
\left[ \begin{array}{c} \phi^{(R)}(x(t),\by(t)) \\ \psi^{(R)}(x(t),\by(t)) \end{array} \right] \;, & x(t) > 0
\end{array} \right\} \,dt + \sqrt{\ee} D \,d\bW(t) \;, \qquad
\left[ \begin{array}{c} x(0) \\ \by(0) \end{array} \right] =
\left[ \begin{array}{c} 0 \\ \by_0 \end{array} \right] \;,
\label{eq:sde2}
\end{equation}
where $\bW(t)$ is an $N$-dimensional standard Brownian motion,
$\ee > 0$ governs the overall strength of the noise,
and the $N \times N$ matrix, $D$, specifies the relative strength of the noise in different directions.
For simplicity we have set $x(0) = 0$, so that the solution is initially located on the switching manifold.
We assume, at least locally, that the functions $\phi^{(L)}$, $\phi^{(R)}$, $\psi^{(L)}$, and $\psi^{(R)}$
are $C^2$ on the closure of their respective half spaces.
The purpose of this assumption is to allow us to study their Taylor series centred at $x=0$.
We write
\begin{equation}
\begin{split}
\phi^{(L)}(x,\by) &= a_L(\by) + c_L(\by) x + O(x^2) \;, \\
\psi^{(L)}(x,\by) &= b_L(\by) + d_L(\by) x + O(x^2) \;, \\
\phi^{(R)}(x,\by) &= -a_R(\by) + c_R(\by) x + O(x^2) \;, \\
\psi^{(R)}(x,\by) &= b_R(\by) + d_R(\by) x + O(x^2) \;,
\end{split}
\label{eq:Taylor}
\end{equation}
where $a_L$, $a_R$, $b_L$ and $b_R$ are $C^2$ functions of $\by$,
and $c_L$, $c_R$, $d_L$ and $d_R$ are $C^1$ functions of $\by$.
Note that $a_L$, $a_R$, $c_L$ and $c_R$ are scalar functions,
whereas $b_L$, $b_R$, $d_L$ and $d_R$ take values in $\mathbb{R}^{N-1}$.

In view of the minus sign in (\ref{eq:Taylor}), subsets of $x=0$ for which
\begin{equation}
a_L(\by) > 0 \;, \qquad
a_R(\by) > 0 \;,
\label{eq:stableSlidingRegion}
\end{equation}
are known as {\em stable sliding regions}, because when $\ee = 0$,
(\ref{eq:sde2}) is a vector field that points toward these regions from both sides.
When $\ee = 0$, forward evolution of a point on a stable sliding region is not defined in the classical sense.
For this reason, as is usual, we employ Filippov's method to define a {\em sliding vector field} on $x=0$,
$\dot{\by}_S = \Omega(\by_S)$,
as the unique convex combination of the limiting left and right vector fields that is tangent to $x=0$ \cite{Fi88,DiBu08,LeNi04}.
Specifically, for (\ref{eq:sde2}) with $\ee = 0$, on $x=0$ we write
\begin{equation}
\left[ \begin{array}{c} \dot{x}_S \\ \dot{\by}_S \end{array} \right] =
(1-\kappa(\by_S)) \left[ \begin{array}{c} a_L(\by_S) \\ b_L(\by_S) \end{array} \right] +
\kappa(\by_S) \left[ \begin{array}{c} -a_R(\by_S) \\ b_R(\by_S) \end{array} \right] \;,
\end{equation}
where $\kappa$ is given by the requirement $\dot{x}_S \equiv 0$, thus $\kappa = \frac{a_L}{a_L+a_R}$.
The sliding vector field is therefore
\begin{equation}
\dot{\by}_S = \frac{a_L(\by_S) b_R(\by_S) + a_R(\by_S) b_L(\by_S)}
{a_L(\by_S) + a_R(\by_S)} \equiv \Omega(\by_S) \;.
\label{eq:Omega}
\end{equation}
We assume that forward orbits of (\ref{eq:sde2}) with $\ee = 0$ evolve on $x=0$, as governed by (\ref{eq:Omega}),
until (\ref{eq:stableSlidingRegion}) ceases to be true and the orbits escape $x=0$.

When $\ee > 0$, an additional definition is not required to define solutions to (\ref{eq:sde2}).
If $\phi^{(L)}$, $\phi^{(R)}$, $\psi^{(L)}$, and $\psi^{(R)}$ are bounded,
then for any $\ee > 0$, (\ref{eq:sde2}) has a unique strong stochastic solution \cite{Fl11,PrSh98,StVa69,KrRo05}.

Our interest is in the transitional PDF of (\ref{eq:sde2})
(i.e.~the PDF of the point $(x(t),\by(t)) \in \mathbb{R}^N$), call it $q(x,\by,t)$.
We assume that at all times in $[0,t]$,
the deterministic solution lies in a stable sliding region,
Therefore $\by_S(t)$ is determined by (\ref{eq:Omega}),
and $a_L(\by_S(s)) > 0$ and $a_R(\by_S(s)) > 0$ for all $s \in [0,t]$.
Below we derive two new explicit limiting expressions for $q(x,\by,t)$.
First we apply Theorem \ref{th:p} to describe marginals of $q$ in the limit $t \to 0$.
Second we perform an asymptotic expansion to describe $q$ in the limit $\ee \to 0$.
This result applies for times, $t \gg \ee$.

\subsection{Asymptotics for small $t$}
\label{sub:SHORT}

To obtain explicit results for small $t$,
we apply two simplifications to (\ref{eq:sde2}).
First, we expect $x(t)$ and $\by(t)$ to not vary greatly over short time frames,
in which case it is reasonable to approximate the drift terms in (\ref{eq:sde2}) by their values at $t=0$.
This gives the approximation
\begin{equation}
\left[ \begin{array}{c} dx(t) \\ d\by(t) \end{array} \right] =
\left\{ \begin{array}{lc}
\left[ \begin{array}{c} a_L \\ b_L \end{array} \right]
\;, & x(t) < 0 \\
\left[ \begin{array}{c} -a_R \\ b_R \end{array} \right]
\;, & x(t) > 0
\end{array} \right\} \,dt + \sqrt{\ee} D \,d\bW(t) \;,
\label{eq:sde4}
\end{equation}
where $a_L$, $a_R$, $b_L$ and $b_R$ are evaluated at $\by(0) = \by_0$.
Second, we assume that the noise in $x$ is independent to the noise in the $\by$. 
Given these assumptions, in this section we derive explicit expressions for the following marginal densities of $q$:
\begin{equation}
q_{\rm orthogonal}(x,t) = \int_{\mathbb{R}^{N-1}} q(x,\by,t) \,d\by \;, \qquad
q_{\rm parallel}(\by,t) = \int_{\mathbb{R}} q(x,\by,t) \,dx \;.
\label{eq:marginals}
\end{equation}
Difficulties in obtaining more general results are discussed at the end of \S\ref{sec:CONC}.
The nomenclature of (\ref{eq:marginals}) was chosen because
the $x$-axis is orthogonal to the switching manifold,
whereas $\by$ comprises of all directions parallel to the switching manifold.

The crucial benefit obtained by the first simplification
is that the right hand side of (\ref{eq:sde4}) is independent of $\by$.
Hence we can decouple (\ref{eq:sde4}) and first write
\begin{equation}
dx(t) =
\left\{ \begin{array}{lc}
a_L \;, & x<0 \\
-a_R \;, & x>0
\end{array} \right\} \,dt + \sqrt{\ee \alpha} \,dW(t) \;,
\qquad x(0) = 0 \;,
\label{eq:dx2}
\end{equation}
where we have contracted the vector noise term into an equivalent scalar noise term
by letting $\alpha = \left( D D^{\sf T} \right)_{11}$.
The transitional PDF of (\ref{eq:dx2}) is
\begin{equation}
q_{\rm orthogonal}(x,t) = \left\{ \begin{array}{lc}
\frac{2}{\ee^3 \alpha^3} \,\re^{\frac{2 a_L x}{\ee \alpha}} \int_0^\infty \int_0^t
h \left( \frac{t-s}{\ee \alpha}, \frac{b}{\ee \alpha}, a_R \right)
h \left( \frac{s}{\ee \alpha}, \frac{b-x}{\ee \alpha}, a_L \right) \,ds \,db \;, & x \le 0 \\
\frac{2}{\ee^3 \alpha^3} \,\re^{\frac{-2 a_R x}{\ee \alpha}} \int_0^\infty \int_0^t
h \left( \frac{t-s}{\ee \alpha}, \frac{b+x}{\ee \alpha}, a_R \right)
h \left( \frac{s}{\ee \alpha}, \frac{b}{\ee \alpha}, a_L \right) \,ds \,db \;, & x \ge 0
\end{array} \right. \;,
\label{eq:qx}
\end{equation}
where $h$ is the first passage time PDF, (\ref{eq:h}).
For a derivation of (\ref{eq:qx}), refer to \cite{KaSh84,KaSh91}.

Second, integration of (\ref{eq:sde4}) yields
\begin{equation}
\by(t) = \by_0 + b_L t + (b_R-b_L)
\int_0^t \chi_{[0,\infty)}(x(s)) \,ds + \sqrt{\ee} \tilde{D} \bW(t) \;,
\label{eq:y}
\end{equation}
where we use $\tilde{D}$ to denote the lower right $(N-1) \times (N-1)$ block of $D$.
Since the scaling, $x \mapsto \frac{x}{\ee \alpha}$, $t \mapsto \frac{t}{\ee \alpha}$,
transforms (\ref{eq:dx2}) to (\ref{eq:dx}) (with $x_0 = 0$),
here the PDF of $\tilde{\tau} = \int_0^t \chi_{[0,\infty)}(x(s)) \,ds$ is
\begin{equation}
p_{\rm scaled}(\tilde{\tau}) = \frac{1}{\ee \alpha} p \left( \frac{\tau}{\ee \alpha}; \frac{t}{\ee \alpha}; 0, a_L, a_R \right) \;.
\end{equation}
Via straight-forward geometric arguments, it follows that the PDF of
$\tilde{\by}(t) = \by_0 + b_L t + (b_R-b_L) \int_0^t \chi_{[0,\infty)}(x(s)) \,ds$,
is given by\removableFootnote{
To derive this, consider the line,
\begin{equation}
l = \left\{ \by_0 + b_L t + (b_R-b_L) \tilde{\tau} ~\big|~ \tilde{\tau} \in \mathbb{R} \right\} \subset \mathbb{R}^{N-1} \;.
\end{equation}
Notice the PDF of $\tilde{\by}$ is zero whenever $\tilde{\by} \notin l$.
The $\tilde{\tau}$-value for the nearest point on $l$ to $\tilde{\by}$ is given by
\begin{equation}
\tilde{\tau} = \frac{(\tilde{\by} - \by_0 - b_L t)^{\sf T} (b_R-b_L)}{|| b_R - b_L ||^2} \;.
\end{equation}
Therefore the distance of $\tilde{\by}$ to $l$ is given by
\begin{equation}
d = \left| \left| \tilde{\by} - \left( \by_0 + b_L t + (b_R-b_L) \tilde{\tau} \right) \right| \right| \;,
\end{equation}
evaluated at the above value of $\tilde{\tau}$.
}
\begin{eqnarray}
\tilde{q}_{\rm parallel}(\tilde{\by},t) &=&
\delta \left( \left| \left|
\tilde{\by} - \by_0 - b_L t -
(b_R-b_L) \frac{(\tilde{\by} - \by_0 - b_L t)^{\sf T} (b_R-b_L)}{|| b_R - b_L ||^2}
\right| \right| \right) \nonumber \\
&&\times~\frac{1}{|| b_R - b_L ||} p_{\rm scaled} \left(
\frac{(\tilde{\by} - \by_0 - b_L t)^{\sf T} (b_R-b_L)}{|| b_R - b_L ||^2}
\right) \;,
\end{eqnarray}
where $\delta$ is the Dirac-delta function.
Finally, the PDF of $\by(t)$,
is the convolution of $\tilde{q}_{\rm parallel}(\tilde{\by},t)$ and the PDF of $\sqrt{\ee} \tilde{D} \bW(t)$:
\begin{equation}
q_{\rm parallel}(\by,t) = \int_{\mathbb{R}^{N-1}} \tilde{q}_{\rm parallel}(\by - \bw,t)
\frac{1}{\left( 2 \pi \ee t \right)^{\frac{N-1}{2}} \sqrt{\det(\gamma)}}
\,\re^{-\frac{1}{2 \ee t} \bw^{\sf T} \gamma^{-1} \bw} \,d\bw \;,
\label{eq:qy}
\end{equation}
where $\gamma = \tilde{D} \tilde{D}^{\sf T}$ (\ref{eq:alphabetagamma}).

\subsection{Asymptotics for small $\ee$}
\label{sub:LONG}

The Fokker-Planck equation for (\ref{eq:sde2}) is
\begin{equation}
q_t = \left\{ \begin{array}{lc}
-\left( \phi^{(L)} q \right)_x
-\left( \psi^{(L)}_i q \right)_{y_i}
+ \frac{\ee \alpha}{2} q_{xx}
+ \ee \beta_i q_{x y_i}
+ \frac{\ee}{2} \gamma_{ij} q_{y_i,y_j} \;, & x < 0 \\
-\left( \phi^{(R)} q \right)_x
-\left( \psi^{(R)}_i q \right)_{y_i}
+ \frac{\ee \alpha}{2} q_{x x}
+ \ee \beta_i q_{x y_i}
+ \frac{\ee}{2} \gamma_{ij} p_{y_i,y_j} \;, & x > 0
\end{array} \right. \;,
\label{eq:fpe7}
\end{equation}
where
\begin{equation}
D D^{\sf T} = \left[ \begin{array}{c|c}
\alpha & \beta^{\sf T} \\ \hline
\beta & \gamma
\end{array} \right] \;,
\label{eq:alphabetagamma}
\end{equation}
and $\alpha \in \mathbb{R}$, $\beta \in \mathbb{R}^{N-1}$ and $\gamma$ is an $(N-1) \times (N-1)$ matrix.
In (\ref{eq:fpe7}), and throughout this section, we use index notation
to abbreviate summations in $i$ and $j$, which always range from $1$ to $N-1$.
The initial and boundary conditions for $q$ are
\begin{eqnarray}
& q(x,\by,t) \to 0 \;, {\rm ~as~} |x|, ||\by|| \to \infty \;, & \\
& q(x,\by,0) = \delta(x) \delta(\by - \by_0) \;, &
\end{eqnarray}
together with a consistency condition at the switching manifold, $x=0$.
This consistency condition is given by the requirement that no probability is gained or lost at $x=0$:
\begin{equation}
e_1^{\sf T} \left( \lim_{x \to 0^-} J - \lim_{x \to 0^+} J \right) = 0 \;,
\label{eq:cc}
\end{equation}
where $J$ denotes the {\em probability current} of (\ref{eq:sde2}) \cite{Sc10,Ga85}.
In the $x$-direction, $J$ is given by
\begin{equation}
e_1^{\sf T} J = \left\{ \begin{array}{lc}
\phi^{(L)} q - \frac{\ee \alpha}{2} q_x
- \frac{\ee}{2} \sum \beta_i q_{y_i} \;, & x < 0 \\
\phi^{(R)} q - \frac{\ee \alpha}{2} q_x
- \frac{\ee}{2} \sum \beta_i q_{y_i} \;, & x > 0
\end{array} \right. \;.
\label{eq:J2}
\end{equation}

To determine $q$ from the above boundary value problem when $\ee$ is small,
a preliminary analysis indicates that a change of variables corresponding to
tracking the system near the deterministic sliding solution is useful.
Specifically we define
\begin{equation}
X = \frac{x}{\ee} \;, \qquad
\bY = \frac{\by - \by_S(t)}{\sqrt{\ee}} \;,
\label{eq:xycheck}
\end{equation}
where $\by_S(t)$ is the deterministic sliding solution determined by (\ref{eq:Omega}).
The motivation for this scaling is discussed in greater detail in earlier work \cite{SiKu13c}.
Under this change of variables, one can obtain an asymptotic approximation from the resulting
Fokker-Planck equation for $Q(X,\bY,t) = q(x,\by,t)$, through a regular series expansion\removableFootnote{
Do we need to multiply by $\ee^{\frac{N+1}{2}}$ in the definition of $Q$?
}
\begin{equation}
Q = Q^{(0)} + \sqrt{\ee} Q^{(1)} + \ee Q^{(2)} + \cdots \;.
\label{eq:qcheckexpansion}
\end{equation}
Here we summarize the steps leading to the asymptotic approximation to $q$ with details given in Appendix \ref{sec:LONG}.
We note that while the calculations for $X < 0$ and $X > 0$ may be treated separately,
the equations for $X > 0$ are the same as those for $X < 0$ 
but with $-a_R$ in place of $a_L$, and $R$'s in place of the remaining $L$'s.

By (\ref{eq:Taylor}) and (\ref{eq:xycheck}),
for $X < 0$ the leading order component of the Fokker-Planck equation (\ref{eq:fpe7}) is
\begin{equation}
\frac{\alpha}{2} Q^{(0)}_{X X} - a_L(\by_S(t)) Q^{(0)}_X = 0 \;.
\end{equation}
Noting that $Q \to 0$ as $X \to -\infty$, the leading order contribution to $Q$ takes the form
\begin{equation}
Q^{(0)}(X,\bY,t) =
f^{(0)}(\bY,t) \,{\rm e}^{\frac{2 a_L(\by_S(t)) X}{\alpha}} \;,
\label{eq:qcheck0}
\end{equation}
for some function $f^{(0)}$.
By also using the analogous form for $X > 0$:
$Q^{(0)}(X,\bY,t) = \linebreak					
f^{(0)}(\bY,t) \,{\rm e}^{\frac{-2 a_R(\by_S(t)) X}{\alpha}}$,
we find that the consistency condition for $Q$ at $X = 0$, (\ref{eq:cc}),
is automatically satisfied for any function $f^{(0)}(\bY,t)$.
Therefore it is necessary to calculate higher order contributions in order to determine $f^{(0)}$,
and hence completely determine the leading order term, $Q^{(0)}$.

For $X < 0$, by taking $O(\sqrt{\ee})$ terms and $O(\ee)$ terms respectively,
we find that the equations for $Q^{(1)}$ and $Q^{(2)}$ take the form
\begin{eqnarray}
\frac{\alpha}{2} Q^{(1)}_{X X} - a_L(\by_S(t)) Q^{(1)}_X
&=& \Phi_1(\bY,t) {\rm e}^{\frac{2 a_L(\by_S(t)) X}{\alpha}} \label{eq:secondOrder} \\
\frac{\alpha}{2} Q^{(2)}_{X X} - a_L(\by_S(t)) Q^{(2)}_X
&=& \left( \Phi_2(\bY,t) + \Phi_3(\bY,t) X \right)
{\rm e}^{\frac{2 a_L(\by_S(t)) X}{\alpha}} \;, \label{eq:thirdOrder}
\end{eqnarray}
where the $\Phi_i$ depend on $f^{(0)}$.
From (\ref{eq:secondOrder}) (and the analogous equation for $X > 0$)
we find that the consistency condition is automatically satisfied at $O(\sqrt{\ee})$.
It is only by using (\ref{eq:thirdOrder}) to apply the consistency condition at $O(\ee)$
that we obtain an equation for $f^{(0)}$, see Appendix \ref{sec:LONG}.
The substitution
\begin{equation}
g^{(0)} = \frac{a_L+a_R}{a_L a_R} \bigg|_{\by = \by_S(t)} f^{(0)} \;,
\label{eq:g0}
\end{equation}
then produces
\begin{equation}
g^{(0)}_t = -\Omega_{i,j} \left( Y_j g^{(0)} \right)_{Y_i}
+ \frac{1}{2} \left(
\frac{r_{Li} r_{Lj} \alpha}{a_L^2}
- \frac{r_{Lj} \beta_i}{a_L}
- \frac{r_{Li} \beta_j}{a_L}
+ \gamma_{ij} \right) \bigg|_{\by = \by_S(t)} g^{(0)}_{Y_i Y_j} \;,
\label{eq:g0PDE}
\end{equation}
where we have introduced
\begin{equation}
r_L = \frac{a_L(b_L-b_R)}{a_L+a_R} \;, \qquad
r_R = -\frac{a_R(b_L-b_R)}{a_L+a_R} \;.
\label{eq:rLrR}
\end{equation}
Equation (\ref{eq:g0PDE}) is the Fokker-Planck equation of the
stochastic differential equation
\begin{equation}
d\bY_{\rm approx}(t) = (D_\by \Omega)(\by_S(t)) \,\bY_{\rm approx}(t) \,dt
+ M(\by_S(t)) \,d\bW(t) \;, \qquad
\bY_{\rm approx}(0) = 0 \;,
\label{eq:linearDiffusionApprox2}
\end{equation}
where $\Omega$ is given by (\ref{eq:Omega}) and 
\begin{equation}
M(\by) = \left[ -\frac{b_L(\by)-b_R(\by)}{a_L(\by)+a_R(\by)} ~\bigg|~ I \right] D \;,
\label{eq:M}
\end{equation}
where $I$ is the $(N-1) \times (N-1)$ identity matrix.
Equation (\ref{eq:linearDiffusionApprox2}) also arises when stochastic averaging arguments are applied
to (\ref{eq:sde2}) with (\ref{eq:xycheck}), \cite{SiKu13}.
Since (\ref{eq:linearDiffusionApprox2}) is a time-dependent Ornstein-Uhlenbeck process \cite{Sc10,Ga09},
the PDF of $\bY_{\rm approx}(t)$ is a zero-mean Gaussian with covariance matrix:
\begin{equation}
\Theta(t) = \int_0^t
{\rm e}^{\int_s^t (D_\by \Omega)(\by_S(u)) \,du}
M(\by_S(s)) M(\by_S(s))^{\sf T}
{\rm e}^{\int_s^t (D_\by \Omega)(\by_S(u))^{\sf T} \,du} \,ds \;.
\label{eq:Theta}
\end{equation}
Since $Q \to 0$ as $|| \by || \to \infty$, we therefore have
\begin{equation}
g^{(0)}(Y,t) = \frac{K(t)}{\sqrt{\det(\Theta(t))}}
\,{\rm e}^{-\frac{1}{2} Y^{\sf T} \Theta(t)^{-1} Y} \;,
\label{eq:g02}
\end{equation}
where $K(t)$ is determined by the requirement that $q$ is normalized.
The combination of 
(\ref{eq:qcheck0}) (and the analogous expression for $X > 0$),
(\ref{eq:g0}) and (\ref{eq:g02}),
provides the following expression for the transitional PDF of (\ref{eq:sde2}):
\begin{equation}
q(x,\by,t) = \frac{1}{(2 \pi \ee)^{\frac{N-1}{2}} \sqrt{\det(\Theta(t))}}
\,{\rm e}^{-\frac{1}{2} \bY^{\sf T} \Theta(t)^{-1} \bY}
\frac{2 a_L a_R}{\alpha \ee (a_L+a_R)} \left\{ \begin{array}{lc}
{\rm e}^{\frac{2 a_L X}{\alpha}} \;, & X < 0 \\
{\rm e}^{-\frac{2 a_R X}{\alpha}} \;, & X > 0
\end{array} \right\} \Bigg|_{\by = \by_S(t)}
+ O \left( \frac{1}{\ee^{\frac{N}{2}}} \right) \;,
\label{eq:qLong}
\end{equation}
where $X$ and $\bY$ are given by (\ref{eq:xycheck}).

Now we consider the approximation obtained by omitting the error term in (\ref{eq:qLong}).
Since $t$ is implicitly assumed to be independent of $\ee$, this approximation applies at times, $t \gg \ee$.
The approximation is piecewise-exponential in $x$ and Gaussian in $\by$.
Moreover,
$x$ and $\by$ are independent in the sense that the value of $x(t)$
does not give us any information about the value of $\by(t)$, and vice-versa.
The mean value of $\by(t)$ is Filippov's sliding solution, $\by_S(t)$.
The covariance of $\by$ limits to zero as $\ee \to 0$, which is consistent with the result of \cite{BuOu09}
that tells us that in the limit $\ee \to 0$, solutions to (\ref{eq:sde2}) are Filippov solutions.
Finally, deviations of $\by$ from the mean are $O(\sqrt{\ee})$,
whereas deviations of $x$ from the mean are $O(\ee)$.
Therefore, there is significantly more variability in the value of $\by(t)$ than in the value of $x(t)$.

\subsection{An illustration of the results for a simple example}
\label{sub:EXAMPLE}

To illustrate the formulas, (\ref{eq:qy}) and (\ref{eq:qLong}), consider
\begin{equation}
\left[ \begin{array}{c} dx(t) \\ dy(t) \end{array} \right] =
\left\{ \begin{array}{lc}
\left[ \begin{array}{c}
-x(t) + y(t) + 1 \\
-x(t) + 1 
\end{array} \right]
\;, & x(t) < 0 \\
\left[ \begin{array}{c}
-x(t) + y(t) - 3 \\
-x(t) - 2
\end{array} \right]
\;, & x(t) > 0
\end{array} \right\} \,dt + \sqrt{\ee}
\left[ \begin{array}{cc} 1 & 0 \\ 0 & \frac{1}{10} \end{array} \right]
\,d\bW(t) \;,
\left[ \begin{array}{c} x(0) \\ y(0) \end{array} \right] =
\left[ \begin{array}{c} 0 \\ 2 \end{array} \right] \;,
\label{eq:sde5}
\end{equation}
which is a two-dimensional example fitting the general form (\ref{eq:sde2}).
Here, $a_L(y) = 1+y$ and $a_R(y) = 3-y$, thus, by (\ref{eq:stableSlidingRegion}),
the interval $(-1,3)$ of the $y$-axis is a stable sliding region.
From (\ref{eq:Omega}) we find that the deterministic sliding solution is
\begin{equation}
y_S(t) = \frac{5}{3} \,\re^{-\frac{3}{4} t} + \frac{1}{3} \;.
\end{equation}

\begin{figure}[t!]
\begin{center}
\includegraphics[height=16.4cm]{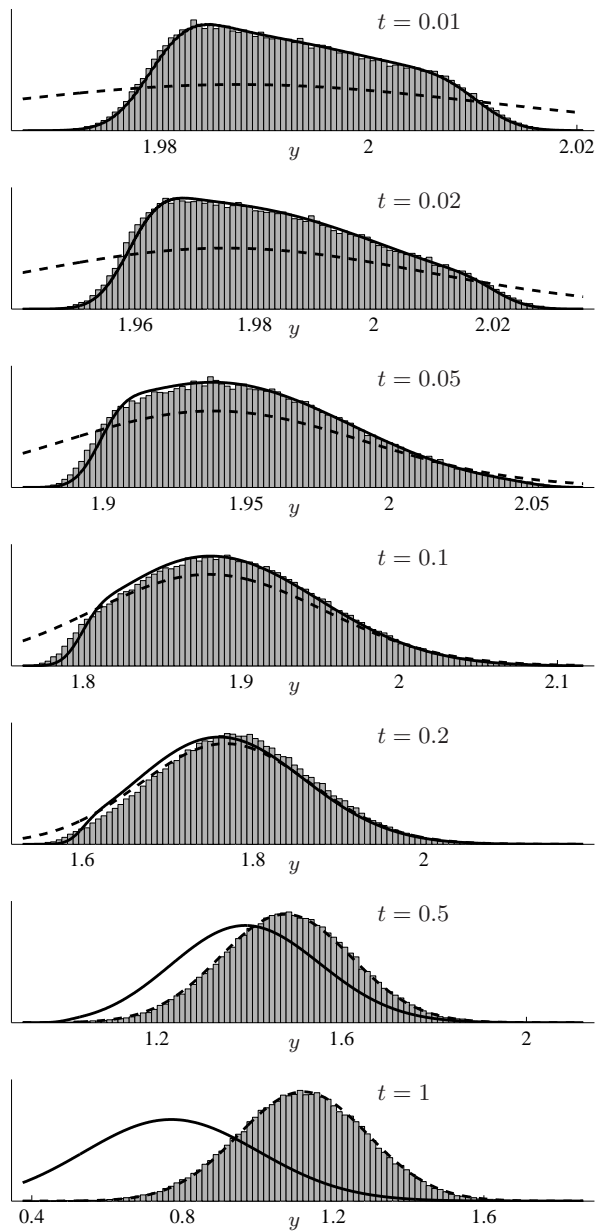}
\caption{
The PDF of $y(t)$ for (\ref{eq:sde5}) with $\ee = 0.1$ at several different values of $t$.
The histograms were computed from numerically computed sample solutions to (\ref{eq:sde5}).
The solid curves show the small time approximation, (\ref{eq:qy}).
The dashed curves show the long time approximation, (\ref{eq:qLong}).
\label{fig:noisySlidingPDFExample}
}
\end{center}
\end{figure}

Fig.~\ref{fig:noisySlidingPDFExample} shows histograms of the value of $y(t)$
as determined from $10^5$ sample solutions to (\ref{eq:sde5})
that were computed numerically using the Euler-Maruyama method of step size $\Delta t = 10^{-5}$.
We have the scaled the histograms to unit area so that they represent the PDF of $y(t)$.
The solid curves show the short time approximation, (\ref{eq:qy}).
The dashed curves show the marginal density in $y$
of the long time approximation, (\ref{eq:qLong}).


As expected, the accuracy of the short time approximation, (\ref{eq:qy}), decreases with increasing $t$,
because the piecewise-constant approximation to (\ref{eq:sde5}) formed by evaluating
the two pieces of the drift at $(x,y) = (0,2)$, as in (\ref{eq:sde4}),
worsens as $(x(t),y(t))$ becomes distant from $(0,2)$.
In contrast, the long time approximation, (\ref{eq:qLong}), is more accurate at larger values of $t$.
This is because (\ref{eq:qLong}) may be thought of as a quasi-steady-state density
that is achieved once sufficient time has passed that sample solutions to (\ref{eq:sde2})
have spent a significant amount of time on each side of $x=0$, with high probability.
In view of the scaling (\ref{eq:xycheck}), this required length of time is of order $\ee$.

\section{Discussion}
\label{sec:CONC}
\setcounter{equation}{0}

This paper gives for the first time the PDF
of the positive occupation time of Brownian motion with two-valued drift.
The PDF, $p$, is given by the complicated expression (\ref{eq:p}).
For large $t$, $p$ is asymptotic to the simpler form (\ref{eq:pLarget}).
In the case that the piecewise drift points toward the origin on both sides, (\ref{eq:pLarget}) is Gaussian.

In \S\ref{sec:SLIDE} we used the result to improve our understanding of the influence
of additive noise on sliding motion.
When the noise amplitude, $\ee$, is zero,
the solution to the general $N$-dimensional system, (\ref{eq:sde2}), slides along $x=0$ for some time.
For $\ee > 0$, we let $q(x,\by,t)$ denote the transitional PDF of (\ref{eq:sde2}).
To approximate $q$ for small $t$,
we constructed a system with piecewise-constant drift, (\ref{eq:sde4}),
by evaluating the drift functions at the initial point.
With this approximation, $x(t)$ is Brownian motion with two-valued drift, (\ref{eq:dx2}),
and $\by(t)$ is a function of the positive occupation time of $x(t)$, (\ref{eq:y}).
In the case that the noise applied to $x$ is independent to the noise applied to $\by$,
it follows immediately from (\ref{eq:y}) that we can write the PDF of $\by(t)$ as
a convolution involving $p$, (\ref{eq:qy}).
The accuracy of the PDF (\ref{eq:qy}) for the general system (\ref{eq:sde2})
correlates with the suitability of approximating the drift terms by their initial values.
For this reason (\ref{eq:qy}) is less accurate at larger times,
and indeed this is evident for the example considered in Fig.~\ref{fig:noisySlidingPDFExample}.

Although the approximate system, (\ref{eq:sde4}), is in general only applicable for small $t$,
it is instructive to apply the long time asymptotic result, (\ref{eq:pLarget}).
For large $t$, the occupation time, $\int_0^t \chi_{[0,\infty)}(x(s)) \,ds$, is asymptotic to 
a Gaussian random variable of mean, $\frac{a_L t}{a_L + a_R}$,
and variance, $\frac{\ee \alpha t}{(a_L+a_R)^2}$.
By substituting this into (\ref{eq:y}), we obtain
\begin{equation}
\mathbb{E}(\by(t)) \sim \by_0 + \Omega(\by_0) t \;, \qquad
{\rm Cov}(\by(t)) \sim \left( \frac{ \alpha (b_R-b_L)(b_R-b_L)^{\sf T}}{(a_L+a_R)^2} + \gamma \right) \ee t \;,
\end{equation}
where $\Omega$ is given by (\ref{eq:Omega}).
Thus the mean of $\by(t)$ is asymptotic to the deterministic sliding solution of the approximate system, (\ref{eq:sde4}),
and the covariance of $\by(t)$ recovers a result of \cite{SiKu13c} that was obtained
directly from the PDF of (\ref{eq:dx}).

For large $t$, the distribution of $\by(t)$ for the original system, (\ref{eq:sde2}),
is approximately Gaussian if the noise amplitude, $\ee$, is small.
Although the approximation, (\ref{eq:qy}), is asymptotically Gaussian for large $t$,
it gives incorrect values for the mean and covariance because the derivation of (\ref{eq:qy}) assumes $t$ is small and,
specifically, does not account for $O(1)$ global changes in the drift as the system evolves.
We derived the small $\ee$ approximation, (\ref{eq:qLong}), which is Gaussian 
with a mean and covariance obtained by appropriately tracking the deterministic sliding solution\removableFootnote{
Technically I have not proved the validity of the asymptotic expansion
which is very difficult, but could probably be achieved by a standard invocation of the maximum principle.
}.
For the example shown in Fig.~\ref{fig:noisySlidingPDFExample},
this approximation shows good agreement to Monte-Carlo simulations
for relatively large values of $t$.

Since $x$ and $\tau$ are strongly correlated,
the product $q_{\rm orthogonal}(x,t) q_{\rm parallel}(\by,t)$ is not a sensible approximation to $q(x,\by,t)$.
It remains to compute the joint PDF of $x$, (\ref{eq:dx}), and $\tau$, (\ref{eq:tau}),
from which we could obtain an accurate expression of $q(x,\by,t)$ for small $t$.
Similarly it remains to compute the joint PDF of $W(t)$ and $\tau$,
which we could use for the general case that the noise applied to $x$
is not independent to the noise applied to $\by$.


\appendix
\section{A calculation of the double inverse Laplace transform}
\label{sec:INVLAP}
\setcounter{equation}{0}

In this section we invert the double Laplace transform (\ref{eq:Xi2}) to obtain $\xi$.
Using the inversion formula, 
$\mathcal{L}^{-1} \left( \frac{1}{\sqrt{s}+k} \right) =
\frac{1}{\sqrt{\pi t}} - k \re^{k^2 t} \,{\rm erfc} \left( k \sqrt{t} \right)$,
together with the standard shifting and scaling rules for Laplace transforms,
the inverse Laplace transform of (\ref{eq:Xi}) with respect to $\omega$ is
\begin{eqnarray}
\int_0^\infty \re^{-\mu \tau}
\xi(\tau,v;a_L,a_R) \,d\tau &=&
\frac{\sqrt{2 \mu + a_R^2} - a_R}{2 \mu} \,\re^{\frac{-a_L^2 v}{2}}
\Bigg( \frac{\sqrt{2}}{\sqrt{\pi v}} -
\left( \sqrt{2 \mu + a_R^2} - a_R - a_L \right) \nonumber \\
&&\times~
\re^{\left( \sqrt{2 \mu + a_R^2} - a_R - a_L \right)^2 \frac{v}{2}}
\,{\rm erfc} \left( \left( \sqrt{2 \mu + a_R^2} - a_R - a_L \right)
\frac{\sqrt{v}}{\sqrt{2}} \right) \Bigg) \;.
\label{eq:Psi}
\end{eqnarray}
To perform the daunting inverse Laplace transform of
(\ref{eq:Psi}), we rewrite it as
\begin{eqnarray}
\int_0^\infty \re^{-\mu \tau}
\xi(\tau,v;a_L,a_R) \,d\tau &=& \re^{\frac{-a_L^2 v}{2}}
\bigg( \frac{1}{\sqrt{\pi v}} \,H_1(\mu;a_R) - H_2(\mu,v;a_L,a_R) \nonumber \\
&&+~
\frac{a_L + 2 a_R}{\sqrt{2}} \,H_1(\mu;a_R) H_2(\mu,v;a_L,a_R) \bigg) \;,
\end{eqnarray}
where
\begin{eqnarray}
H_1(\mu;a_R) &=& \frac{\sqrt{2 \mu + a_R^2} - a_R}{\sqrt{2} \mu} \;, \\
H_2(\mu,v;a_L,a_R) &=& 
\re^{\left( \sqrt{2 \mu + a_R^2} - a_R - a_L \right)^2 \frac{v}{2}}
\,{\rm erfc} \left( \left( \sqrt{2 \mu + a_R^2} - a_R - a_L \right)
\frac{\sqrt{v}}{\sqrt{2}} \right) \;,
\end{eqnarray}
so that
\begin{equation}
\xi(\tau,v;a_L,a_R) = \re^{\frac{-a_L^2 v}{2}}
\left( \frac{1}{\sqrt{\pi v}} \,h_1(\tau;a_R) - h_2(\tau,v;a_L,a_R) +
\frac{a_L + 2 a_R}{\sqrt{2}} \,h_1(\tau;a_R) * h_2(\tau,v;a_L,a_R) \right) \;,
\label{eq:psi}
\end{equation}
where $h_1$ and $h_2$ denote the inverse Laplace transforms
of $H_1$ and $H_2$ respectively,
and $*$ denotes convolution with respect to $\tau$.
By writing
$H_1(\mu;a_R) = -\frac{a_R}{\sqrt{2} \mu} + \frac{\sqrt{2}}{\sqrt{2 \mu + a_R^2}} +
\frac{a_R^2}{\sqrt{2} \mu \sqrt{2 \mu + a_R^2}}$,
and using, in particular\removableFootnote{
Also: $\mathcal{L}^{-1} \left( \frac{1}{\sqrt{s+k}} \right) = \frac{\re^{-kt}}{\sqrt{\pi t}}$.
},
$\mathcal{L}^{-1} \left( \frac{1}{s \sqrt{s+k}} \right) =
\frac{1}{\sqrt{k}} \left( 1 - {\rm erfc} \left( \sqrt{kt} \right) \right)$,
we obtain
\begin{equation}
h_1(\tau;a_R) = \frac{\re^{\frac{-a_R^2 \tau}{2}}}{\sqrt{\pi \tau}} -
\frac{a_R}{\sqrt{2}} \,{\rm erfc} \left( \frac{a_R \sqrt{\tau}}{\sqrt{2}} \right) \;.
\label{eq:g1}
\end{equation}
Also by using
$\mathcal{L}^{-1}
\left( \re^{(\sqrt{s}+k)^2} \,{\rm erfc} \left( \sqrt{s} + k \right) \right) =
\frac{1}{2 \pi t^{\frac{3}{2}}}
\int_0^\infty z \re^{\frac{-z^2}{4} \left( 1 + \frac{1}{t} \right) }
\re^{-kz} \,dz$ (obtained by combining formulas in \cite{RoKa66})\removableFootnote{
I created this identity by applying the obscure formula
(given on page 171 of \cite{RoKa66})
\begin{equation}
\mathcal{L}^{-1} \left( F \left( \sqrt{s} \right) \right) =
\frac{1}{2 \sqrt{\pi} t^{\frac{3}{2}}}
\int_0^\infty z \re^{\frac{-z^2}{4t}} f(z) \,dz \;,
\end{equation}
(where $F$ is an arbitrary function and $f(t) = \mathcal{L}^{-1} \left( F(s) \right)$)
to
\begin{equation}
\mathcal{L}^{-1} \left( \re^{(s+k)^2} \,{\rm erfc} (s+k) \right) =
\frac{1}{\sqrt{\pi}} \re^{-kt} \re^{\frac{-t^2}{4}} \;.
\end{equation}
},
we obtain
\begin{equation}
h_2(\tau,v;a_L,a_R) =
\frac{\re^{\frac{-a_R^2 \tau}{2}}}{\sqrt{\pi} (\tau+v)}
\left( \frac{\sqrt{v}}{\sqrt{\pi \tau}} +
\frac{(a_L+a_R)v}{\sqrt{2(\tau+v)}}
\re^{\frac{(a_L+a_R)^2 \tau v}{2(\tau+v)}}
\,{\rm erfc} \left( \frac{-(a_L+a_R) \sqrt{\tau v}}{\sqrt{2(\tau+v)}} \right) \right) \;.
\label{eq:g2}
\end{equation}
After an involved algebraic simplification\removableFootnote{
Note that
\begin{equation}
h_1(\tau;a_R) = \frac{2}{a_R^2} h_3'(\tau;a_R) + h_3(\tau;a_R) \;,
\end{equation}
where
\begin{equation}
h_3(\tau;a_R) = -\frac{a_R}{\sqrt{2}} \,{\rm erfc} \left( \frac{a_R \sqrt{\tau}}{\sqrt{2}} \right) \;.
\end{equation}
To evaluate the convolution appearing in (\ref{eq:psi}) we first find that
\begin{eqnarray}
\frac{2}{a_R^2} h_3'(\tau;a_R) * h_2(\tau,v;a_L,a_R) &=&
\frac{\sqrt{v} \,\re^{\frac{-a_R^2 \tau}{2}}}{\pi^{\frac{3}{2}}}
\int_0^\tau \frac{1}{\sqrt{z(\tau-z)} (z+v)} \,dz \nonumber \\
&&+~
\frac{(a_L+a_R) v \,\re^{\frac{-a_R^2 \tau}{2}}}{\sqrt{2} \pi}
\int_0^\tau \frac{1}{\sqrt{2(\tau-z)} (z+v)^{\frac{3}{2}}}
\,\re^{\frac{(a_L+a_R)^2 z v}{2(z+v)}}
\,{\rm erfc} \left( \frac{-(a_L+a_R) \sqrt{z v}}{\sqrt{2(z+v)}} \right) \,dz \nonumber \\
&=&
\frac{\re^{\frac{-a_R^2 \tau}{2}}}{\sqrt{\pi} \sqrt{\tau+v}}
\,\re^{\frac{a_L+a_R)^2 \tau v}{2(\tau+v)}}
\,{\rm erfc} \left( \frac{-(a_L+a_R) \sqrt{\tau v}}{\sqrt{2(\tau+v)}} \right) \;,
\end{eqnarray}
where the first integral can evaluated by using the
convenient exact expression for the antiderivative:
$\int \frac{1}{\sqrt{z(\tau-z)} (z+v)} \,dz =
\frac{2}{\sqrt{v(\tau+v)}}
\tan^{-1} \left( \frac{\sqrt{(\tau+v)z}}{\sqrt{v} \sqrt{\tau-z}} \right) + C$,
and the second integral succumbs to the substitution $x = \frac{zv}{z+v}$
and further appropriate manipulation.

We then use the following nifty identity (true for arbitrary functions)
\begin{equation}
h_3(\tau;a_R) * h_2(\tau,v;a_L,a_R) =
\int_0^\tau h_3'(z;a_R) * h_2(z,v;a_L,a_R) + h_3(0;a_R) h_2(z,v;a_L,a_R) \,dz \;,
\end{equation}
to arrive at
\begin{equation}
\frac{a_L + 2 a_R}{\sqrt{2}} \,\re^{\frac{-a_L^2 v}{2}} 
h_3(\tau;a_R) * h_2(\tau,v;a_L,a_R) = \mathcal{F}(\tau;t;a_L,a_R) \;.
\end{equation}
}
that we omit for brevity,
the combination of (\ref{eq:p2}), (\ref{eq:psi}), (\ref{eq:g1}) and (\ref{eq:g2})
yields (\ref{eq:p}) with (\ref{eq:Fcal}) as required.

\section{Details of the asymptotic expansion}
\label{sec:LONG}
\setcounter{equation}{0}




For $X < 0$, the Fokker-Planck equation for $Q(X,\bY,t)$ is\removableFootnote{
In particular we have used
\begin{eqnarray}
\phi^{(L)}(\by) &=&
a_L(\by) + c_L(\by) x + O \left( x^2 \right) \nonumber \\
&=& a_L \left( \by_S(t) + \sqrt{\ee} \by \right)
+ c_L \left( \by_S(t) + \sqrt{\ee} \by \right) \ee X + O(\ee^2) \nonumber \\
&=& a_L + \sqrt{\ee} \frac{\partial a_L}{\partial y_i} Y_i
+ \frac{\ee}{2} \frac{\partial^2 a_L}{\partial y_i \partial y_j}
Y_i Y_j
+ \ee c_L X + O \left( \ee^{\frac{3}{2}} \right) \;,
\end{eqnarray}
where in the last line the coefficients are evaluated at $\by = \by_S(t)$.
Also
\begin{equation}
q_t = Q_X \frac{\partial X}{\partial t}
+ Q_{Y_i} \frac{\partial Y_i}{\partial t}
+ Q_t
= Q_t - \frac{1}{\sqrt{\ee}} Q_{Y_i}
\frac{d \by_{{\rm sl},i}}{dt}
= Q_t - \frac{1}{\sqrt{\ee}} \Omega_i Q_{Y_i} \;.
\end{equation}
}
\begin{eqnarray}
Q_t &=&
-\frac{1}{\ee} a_L Q_X
- \frac{1}{\sqrt{\ee}} a_{L,i} Y_i Q_X
- \frac{1}{2} a_{L,ij}
Y_i Y_j Q_X
- c_L \left( X Q \right)_X
- \frac{1}{\sqrt{\ee}} b_{Li} Q_{Y_i}
- b_{Li,j} \left( Y_j Q \right)_{Y_i} \nonumber \\
&&+~\frac{\alpha}{2 \ee} Q_{X X}
+ \frac{1}{\sqrt{\ee}} \beta_i Q_{X Y_i}
+ \frac{1}{2} \gamma_{ij} Q_{Y_i Y_j}
+ \frac{1}{\sqrt{\ee}} \Omega_i Q_{Y_i} + O(\sqrt{\ee}) \;,
\label{eq:fpe8}
\end{eqnarray}
where for the remainder of this section, unless otherwise specified, the coefficients, $a_J, \ldots, d_J$,
and their derivatives, are evaluated at $\by = \by_S(t)$,
and we have abbreviated derivatives with respect to components of $\by$
with indices in subscripts following commas:
\begin{equation}
a_{L,i} \equiv \frac{\partial a_L}{\partial y_i} \;, \qquad
a_{L,ij} \equiv \frac{\partial^2 a_L}{\partial y_i \partial y_j} \;, \qquad
b_{Li,j} \equiv \frac{\partial b_{Li}}{\partial y_j} \;.
\end{equation}
By multiplying (\ref{eq:fpe8}) (and the analogous expression for $X > 0$)
through by $\ee$ and collecting terms of the same order, we arrive at
\begin{equation}
\frac{\alpha}{2} Q_{X X} - a_L Q_X =
\left\{ \begin{array}{lc}
\sqrt{\ee} \left( a_{L,i} Y_i Q_X
+ r_{Li} Q_{Y_i}
- \beta_i Q_{X Y_i} \right) + \ee \Big(
\frac{1}{2} a_{L,ij} Y_i Y_j Q_X
+ c_L \left( X Q \right)_X \\
\qquad+~b_{Li,j} \left( Y_j Q \right)_{Y_i}
- \frac{1}{2} \gamma_{ij} Q_{Y_i Y_j}
+ Q_t \Big) + O \left( \ee^{\frac{3}{2}} \right) \;, & X < 0 \\
\sqrt{\ee} \left( -a_{R,i} Y_i Q_X
+ r_{Ri} Q_{Y_i}
- \beta_i Q_{X Y_i} \right) + \ee \Big(
-\frac{1}{2} a_{R,ij} Y_i Y_j Q_X
+ c_R \left( X Q \right)_X \\
\qquad+~b_{Ri,j} \left( Y_j Q \right)_{Y_i}
- \frac{1}{2} \gamma_{ij} Q_{Y_i Y_j}
+ Q_t \Big) + O \left( \ee^{\frac{3}{2}} \right) \;, & X > 0
\end{array} \right.
\label{eq:fpe9}
\end{equation}
where we have substituted (\ref{eq:rLrR}).
In terms of the scaled variables, the consistency condition, (\ref{eq:cc}), may be written as
\begin{eqnarray}
\lim_{X \to 0^-} Q_X(X,\bY,t) -
\lim_{X \to 0^+} Q_X(X,\bY,t) &=&
\frac{2}{\alpha} \bigg( a_L + a_R
+ \sqrt{\ee} \left( a_{L,i} + a_{R,i} \right) Y_i \nonumber \\
&&+~\frac{\ee}{2} \left( a_{L,ij} + a_{R,ij} \right) Y_i Y_j
+ O \left( \ee^{\frac{3}{2}} \right)
\bigg) Q(0,\bY,t) \;.
\label{eq:cc3}
\end{eqnarray}

The leading order term, $Q^{(0)}$, of the asymptotic expansion of $Q$ is given by (\ref{eq:qcheck0}).
The equation for the $O(\sqrt{\ee})$ correction to $Q^{(0)}$ is given by
\begin{equation}
\frac{\alpha}{2} Q^{(1)}_{X X} - a_L Q^{(1)}_X = \Phi_1(\bY,t) {\rm e}^{\frac{2 a_L X}{\alpha}} \;,
\end{equation}
where
\begin{equation}
\Phi_1(\bY,t) = \frac{2 a_L}{\alpha} a_{L,i} Y_i f^{(0)}
+ r_{Li} f^{(0)}_{Y_i}
- \frac{2 a_L}{\alpha} \beta_i f^{(0)}_{Y_i} \;.
\end{equation}
Consequently,
\begin{equation}
Q^{(1)}(X,\bY,t) =
\left( f^{(1)}(\bY,t) + \frac{1}{a_L} \Phi_1(\bY,t) X \right)
{\rm e}^{\frac{2 a_L X}{\alpha}} \;,
\label{eq:qcheck1}
\end{equation}
for some function $f^{(1)}$.
Since the consistency condition (\ref{eq:cc3}) is automatically satisfied to $O(1)$ and $O(\sqrt{\ee})$,
$f^{(0)}$ and $f^{(1)}$ are undetermined.
The equation for $Q^{(2)}$ is
\begin{equation}
\frac{\alpha}{2} Q^{(2)}_{X X} - a_L Q^{(2)}_X =
\left( \Phi_2(\bY,t) + \Phi_3(\bY,t) X \right) {\rm e}^{\frac{2 a_L X}{\alpha}} \;,
\end{equation}
where
\begin{eqnarray}
\Phi_2(\bY,t) &=&
\frac{2 a_L}{\alpha} a_{L,i} Y_i f^{(1)}
+ \frac{2}{\alpha} a_{L,i} a_{L,j} Y_i Y_j f^{(0)}
+ \frac{1}{a_L} a_{L,i} r_{Lj} Y_i f^{(0)}_{Y_j} \nonumber \\
&&-~\frac{2}{\alpha} a_{L,i} \beta_j Y_i f^{(0)}_{Y_j}
+ r_{Li} f^{(1)}_{Y_i}
- \frac{2 a_L}{\alpha} \beta_i f^{(1)}_{Y_i}
- \frac{2}{\alpha} a_{L,j} \beta_i \left( Y_j f^{(0)} \right)_{Y_i} \nonumber \\
&&-~\frac{1}{a_L} r_{Lj} \beta_i f^{(0)}_{Y_i Y_j}
+ \frac{2}{\alpha} \beta_i \beta_j f^{(0)}_{Y_i Y_j}
+ \frac{a_L}{\alpha} a_{L,ij} Y_i Y_j f^{(0)}
+ c_L f^{(0)} \nonumber \\
&&+~b_{Li,j} \left( Y_j f^{(0)} \right)_{Y_i}
- \frac{1}{2} \gamma_{ij} f^{(0)}_{Y_i Y_j}
+ f^{(0)}_t \;, \\
\Phi_3(\bY,t) &=&
\frac{4 a_L}{\alpha^2} a_{L,i} a_{L,j} Y_i Y_j f^{(0)}
+ \frac{2}{\alpha} a_{L,i} r_{Lj} Y_i f^{(0)}_{Y_j}
- \frac{4 a_L}{\alpha^2} a_{L,i} \beta_j Y_i f^{(0)}_{Y_j} \nonumber \\
&&+~\frac{2}{\alpha} a_{L,j} r_{Li} \left( Y_j f^{(0)} \right)_{Y_i}
+ \frac{1}{a_L} r_{Li} r_{Lj} f^{(0)}_{Y_i Y_j}
- \frac{2}{\alpha} r_{Li} \beta_j f^{(0)}_{Y_i Y_j} \nonumber \\
&&-~\frac{4 a_L}{\alpha^2} a_{L,j} \beta_i \left( Y_j f^{(0)} \right)_{Y_i}
- \frac{2}{\alpha} r_{Li} \beta_j f^{(0)}_{Y_i Y_j}
+ \frac{4 a_L}{\alpha^2} \beta_i \beta_j f^{(0)}_{Y_i Y_j} \nonumber \\
&&+~\frac{2 a_L c_L}{\alpha} f^{(0)}
+ \frac{2}{\alpha} a_{L,i} \Omega_i f^{(0)} \;.
\end{eqnarray}
Consequently
\begin{equation}
Q^{(2)}(X,\bY,t) = \left( f^{(2)}(\bY,t)
+ \left( \frac{1}{a_L} \Phi_2(\bY,t) - \frac{\alpha}{2 a_L^2} \Phi_3(\bY,t) \right) X
+ \frac{1}{2 a_L} \Phi_3(\bY,t) X^2 \right) 
{\rm e}^{\frac{2 a_L X}{\alpha}} \;,
\label{eq:qcheck2}
\end{equation}
for some function $f^{(2)}$.
By applying the consistency condition at $O(\ee)$,
after an involved algebraic reduction\removableFootnote{
Regarding (\ref{eq:cc3}),
\begin{eqnarray}
{\rm RHS} &=&
\frac{2}{\alpha} (a_L + a_R) f^{(0)}
+ \frac{2 \sqrt{\ee}}{\alpha} \left( (a_L+a_R) f^{(1)}
+ (a_{L,i}+a_{R,i}) Y_i f^{(0)} \right) \nonumber \\
&&+~\frac{2 \ee}{\alpha} \left( (a_L+a_R) f^{(2)}
+ (a_{L,i}+a_{R,i}) Y_i f^{(1)}
+ \frac{1}{2} (a_{L,ij}+a_{R,ij}) Y_i Y_j f^{(0)} \right)
+ O \left( \ee^{\frac{3}{2}} \right) \;.
\end{eqnarray} 
By combining with the analogous solution for $X \ge 0$,
and repeatedly using identity
\begin{equation}
\frac{r_L}{a_L} + \frac{r_R}{a_R} = 0 \;,
\end{equation}
regarding (\ref{eq:cc3}),
\begin{eqnarray}
{\rm LHS} &=& {\rm RHS} + \ee \Bigg(
\left( -\frac{a_{L,j} r_{Li}}{a_L^2} - \frac{a_{R,j} r_{Ri}}{a_R^2}
+ \frac{b_{Li,j}}{a_L} + \frac{b_{Ri,j}}{a_R} \right)
\left( Y_j f^{(0)} \right)_{Y_i} \nonumber \\
&&-~\left( \left(
\frac{r_{Li} r_{Lj}}{a_L^3}
+ \frac{r_{Ri} r_{Rj}}{a_R^3} \right) \frac{\alpha}{2}
+ \left( -\frac{r_{Lj}}{a_L^2} + \frac{r_{Rj}}{a_R^2} \right) \beta_i
+ \left( \frac{1}{a_L} + \frac{1}{a_R} \right) \frac{\gamma_{ij}}{2} \right)
f^{(0)}_{Y_i Y_j} \nonumber \\
&&+~\left( \frac{1}{a_L} + \frac{1}{a_R} \right) f^{(0)}_t
- \left( \frac{a_{L,i} \Omega_i}{a_L^2} + \frac{a_{R,i} \Omega_i}{a_R^2} \right) f^{(0)}
\Bigg)
+ O \left( \ee^{\frac{3}{2}} \right) \;.
\label{eq:ccConstruct}
\end{eqnarray}
To reduce the first line of (\ref{eq:ccConstruct}) we note:
\begin{eqnarray}
&& -\frac{a_{L,j} r_{Li}}{a_L^2} - \frac{a_{R,j} r_{Ri}}{a_R^2}
+ \frac{b_{Li,j}}{a_L} + \frac{b_{Ri,j}}{a_R} =
\frac{a_L b_{Ri,j} + a_{L,j} r_{Ri} + a_R b_{Li,j} + a_{R,j} r_{Li}}
{a_L a_R} \nonumber \\
&=& \frac{a_L+a_R}{a_L a_R} \left(
\frac{a_L b_{Ri,j} + a_{L,j} b_{Ri} + a_R b_{Li,j} + a_{R,j} b_{Li}}{a_L+a_R}
- \frac{(a_L b_{Ri} + a_R b_{Li})(a_{L,j}+a_{R,j})}{(a_L+a_R)^2} \right) \\
&=& \frac{a_L+a_R}{a_L a_R} \Omega_{i,j} \;.
\end{eqnarray}
The second line of (\ref{eq:ccConstruct}) reduces simply:
\begin{equation}
\left( \frac{r_{Li} r_{Lj}}{a_L^3}
+ \frac{r_{Ri} r_{Rj}}{a_R^3} \right) \frac{\alpha}{2}
+ \left( -\frac{r_{Lj}}{a_L^2} + \frac{r_{Rj}}{a_R^2} \right) \beta_i
+ \left( \frac{1}{a_L} + \frac{1}{a_R} \right) \frac{\gamma_{ij}}{2}
= \frac{a_L+a_R}{2 a_L a_R} \left(
\frac{r_{Li} r_{Lj} \alpha}{a_L^2}
- \frac{2 r_{Lj} \beta_i}{a_L} + \gamma_{ij} \right) \;.
\end{equation}
Finally the third line simplifies conveniently:
\begin{equation}
\left( \frac{1}{a_L} + \frac{1}{a_R} \right) f^{(0)}_t
- \left( \frac{a_{L,i} \Omega_i}{a_L^2} + \frac{a_{R,i} \Omega_i}{a_R^2} \right) f^{(0)}
= \left( \frac{a_L+a_R}{a_L a_R} f^{(0)} \right)_t \;.
\end{equation}
Therefore, (\ref{eq:ccConstruct}) reduces to
\begin{eqnarray}
{\rm LHS} &=& {\rm RHS} + \ee \Bigg(
\frac{a_L+a_R}{a_L a_R} \Omega_{i,j} \left( Y_j f^{(0)} \right)_{Y_i} \nonumber \\
&&-~\frac{a_L+a_R}{2 a_L a_R} \left(
\frac{r_{Li} r_{Lj} \alpha}{a_L^2}
- \frac{2 r_{Lj} \beta_i}{a_L} + \gamma_{ij} \right) f^{(0)}_{Y_i Y_j}
+ \left( \frac{a_L+a_R}{a_L a_R} f^{(0)} \right)_t \Bigg)
+ O \left( \ee^{\frac{3}{2}} \right) \;.
\end{eqnarray}
}
we arrive at (\ref{eq:g0PDE}).
The observation that (\ref{eq:g0PDE}) is the Fokker-Planck equation
for the stochastic differential equation, (\ref{eq:linearDiffusionApprox2}),
leads to the desired result (\ref{eq:qLong}) as described in the main text.

\end{document}